\documentclass[11pt,draft]{article}

\usepackage{amsfonts,amsmath,amssymb,amsthm}

\newcommand{\rset}{\mathbf{R}}

\newcommand{\nset}{\mathbf{N}}

\newcommand{\ind}{\mathbf{1}}

\newcommand{\e}{\mathbb{E}}

\newcommand{\m}{\mathcal}
\usepackage{comment}
\usepackage[a4paper,body={160mm,230mm},centering]{geometry}
\newtheorem{thm}{Theorem}[section]
\newtheorem{lemme}[thm]{Lemma}
\newtheorem{prop}[thm]{Proposition}
\newtheorem{cor}[thm]{Corollary}
\theoremstyle{definition}
\newtheorem{df}[thm]{Definition}
\newtheorem{hyp}{Assumption}

\theoremstyle{remark}
\newtheorem{rem}[thm]{Remark}
\def\ds{\begin{displaystyle}}
\def\eds{\end{displaystyle}}
\def\dis{\displaystyle }
\def\<{\langle }
\def\>{\rangle }
\newcommand{\nl}{\nolimits}

\def\R{\mathbf{ R}}
\def\N{\mathbf{N}}

\def\E{\mathbb E}

\def\P{\mathbb P}
\def\Q{\mathbb Q}

\def\calb{{\cal B}}

\def\calf{{\cal F}}
\def\calg{{\cal G}}

\def\calp{{\cal P}}

\def\call{{\cal L}}

\def\by{{\overline{Y}}}
\def\bz{{\overline{Z}}}

\def\bW{{\overline{W}}}
\def\wy{{\widehat{Y}}}
\def\wz{{\widehat{Z}}}

\title{\bf Quadratic BSDEs with random terminal time and elliptic PDEs in infinite dimension.}

\author{Philippe Briand\\[.3em]
\normalsize\it
IRMAR, Universit\'{e} Rennes 1, 35042 Rennes Cedex, FRANCE\\
\normalsize\tt philippe.briand@univ-rennes1.fr \and
Fulvia Confortola\\[.3em]
\normalsize\it
Dipartimento di Matematica e Applicazioni,\\
\normalsize\it Universit\`a di Milano-Bicocca\\
\normalsize\it Via R. Cozzi 53 - Edificio U5 - 20125 Milano, Italy \\
\normalsize\tt fulvia.confortola@unimib.it}
\date{}

\begin{document}

\maketitle


\bigskip

\centerline{\bf Abstract} In this paper we study one dimensional
backward stochastic differential equations (BSDEs) with random
terminal time not necessarily bounded or finite when the generator
$F(t,Y,Z)$ has a quadratic growth in $Z$. We provide existence and
uniqueness of a bounded solution of such BSDEs and, in the case of
infinite horizon, regular dependence on parameters. The obtained
results are then applied to prove existence and uniqueness of a
mild solution to elliptic partial differential equations in
Hilbert spaces. Finally we show an application to a control
problem.

\section{Introduction}

Let $\tau$ be a stopping time which is not necessarily bounded or
finite. We look for a pair of processes $(Y_t,Z_t)_{t \geq 0}$
progressively measurable which satisfy $\forall t \geq 0, \forall
T \geq t$
\begin{equation}\label{bsde-stime-intro}\left\{
\begin{array}{l}
 \displaystyle    Y_{t \wedge \tau}= Y_{T\wedge \tau} +\int_{t\wedge \tau}^{T\wedge \tau} F(s,Y_s,Z_s) - \int_{t\wedge \tau}^{T\wedge \tau} Z_s dW_s\\
\\
\displaystyle Y_{\tau}= \xi \mbox{ on }
\{\tau<\infty\}\end{array}\right.
\end{equation}
where $W$ is a cylindrical Wiener process in some infinite
dimensional Hilbert space $\Xi$ and the generator $F$ has
quadratic growth with respect to the variable $z$. Moreover the
terminal condition $\xi$ is $\m F_{\tau}$-measurable and bounded.
We limit ourselves to the case in which $(Y_t)_{t \geq 0}$ is
one-dimensional and we look for a solution $(Y_t,Z_t)_{t \geq 0}$
such that $(Y_t)_{t \geq 0}$ is a bounded process and $(Z_t)_{t
\geq 0}$ is a process with values in the space of the
Hilbert-Schmidt operator from $\Xi$ to $\rset$ such that
$\E\left(\int_0^{t \wedge \tau} |Z_s|^2ds\right)< \infty, \forall
t \geq 0$.

BSDEs with random terminal time have been treated by several
authors (see for instance \cite{Pa99}, \cite{DaPa97}, \cite{BH98},
\cite{Ro}) when the generator is Lipschitz, or monotone and with
suitable growth with respect to $y$, but Lipschitz with respect to
$z$. Kobylanski \cite{Kob00} deals with a real BSDE with quadratic
generator with respect to $z$ and with random terminal time. She
requires that the stopping time is bounded or $\P$-a.s finite.
We generalize in a certain sense the result of Kobylanski, but, to
obtain the existence and uniqueness of the solution to
(\ref{bsde-stime-intro}) for a general stopping time, we have to
require stronger assumption on the generator. In particular it has
to be strictly monotone with respect to $y$.

We follow the techniques introduced by Briand and Hu in
\cite{BH98}, and used successively by Royer \cite{Ro}, based upon
an approximation procedure and on Girsanov transform. We can use
this strategy even if, under our assumptions, the generator is not
Lipschitz with respect to $z$. The main idea is to exploit the
theory of BMO-martingales. It is indeed known that if $(Y,Z)$
solves a quadratic BSDE with bounded (or $\P$-a.s.) finite final
time then $\displaystyle \int_0^{\cdot} Z_s\, dW_s$ is a
BMO--martingale (see \cite{HIM05}).

Then the result on BSDE is exploited to study existence and
uniqueness of a mild solution (see Section \ref{sec-Kolm} for the
definition) to the following elliptic partial differential
equation in Hilbert space $H$
\begin{equation}\label{pdeintroduction}
{\cal L}u(x)+F(x,u(x),\nabla u(x)\sigma)=0,\quad x\in H,
\end{equation}
where $F$ is a function from $H\times \mathbb R\times \Xi^*$ to
$\mathbb R$ strictly monotone with respect the second variable and
with quadratic growth in the gradient of the solution and ${\cal
L}$ is the second order operator:
$${\cal L}\phi(x)=\frac{1}{2} Trace(\sigma \sigma^*\nabla^2\phi(x))+\langle
Ax,\nabla\phi(x)\rangle+\langle b(x),\nabla\phi(x)\rangle.$$ $H$
is an Hilbert space, $A$ is the generator of a strongly continuous
semigroup of bounded linear operators $(e^{tA})_{t\ge 0}$ in $H$,
$b$ is a function with values in $H$ and $\sigma$ belongs to
$L(\Xi,H)$- the space of linear bounded operator from $\Xi$ to $K$
satisfying appropriate Lipschitz
 conditions.

Existence and uniqueness of a mild solution of equation
(\ref{pdeintroduction}) in infinite dimensional spaces have been
recently studied by several authors employing different techniques
(see \cite{Ce}, \cite{GoRo}, \cite{dz3} and \cite{fede}).

In \cite{FT3} (following several papers
   dealing with finite dimensional situations, see, for instance \cite{BuPe},
\cite{DaPa97} and \cite{Pa98}) the solution of equation
(\ref{pdeintroduction}) is represented using a Markovian
forward-backward system of equations
\begin{equation}\label{fbsde1}\left\{
\begin{array}{l}
\displaystyle dX_s=AX_sds+b(X_s)ds+\sigma(X_s)dW_s,\quad s\geq 0 \\
\displaystyle dY_s= - F(X_s,Y_s,Z_s)ds+Z_s dW_s,\quad s\geq 0\\
X_0=x
\end{array}\right.
\end{equation}
where $F$ is Lipschitz with respect to $y$ and $z$ and monotone in
$y$, but with monotonicity constant large. A such limitation has
then been removed under certain conditions in \cite{HT}, still
assuming $F$ Lipschitz with respect to $z$, strictly monotone and
with arbitrary growth with respect to $y$. We follow the same
approach to deal with mild solution to
 (\ref{pdeintroduction}) when the coefficient $F$ is strictly monotone in the second variable (there are not conditions
 on its monotonicity constant) and has quadratic growth in the gradient of the solution.
The main technical point here will be proving differentiability of
the bounded solution of the
 backward equation in system (\ref{fbsde1}) with respect to the initial datum $x$ of the
 forward equation. To obtain this result we follow \cite{HT}. The proof is based on an a-priori bound for suitable
    approximations of the equations for the gradient of $Y$ with respect to
    $x$. We use again classical result on BMO-martingales.

In the last part of the paper we apply the above result to an
optimal control problem with state equation:

\begin{equation}\label{sdecontrol}\left\{\begin{array}{l}
d X_{\tau}=AX_{\tau}d \tau+b(X_\tau )d \tau+ \sigma
r(X_{\tau},u_{\tau}) d\tau+ \sigma dW_{\tau},\\
X_0=x\in H,
\end{array}\right.
\end{equation}
where $u$ denotes the control process, taking values in a given
closed subset ${\cal U}$ of a Banach space $U$. The control
problem consists of minimizing an infinite horizon cost functional
of the form
$$J(x,u)=\mathbb E\int_0^\infty e^{-\lambda\sigma}g(X_\sigma^u,u_\sigma)d\sigma.$$
We suppose that $r$ is a function with values in $\Xi^*$ with
linear growth in $u$ and $g$ is a given real function with
quadratic growth in $u$. $\lambda$ is any positive number. We
assume that neither $\m U$ nor $r$ is bounded: in this way the
Hamiltonian corresponding to the control problem has quadratic
growth in the gradient of the solution and consequently the
associated BSDE has quadratic growth in the variable $Z$. The
results obtained on equation (\ref{pdeintroduction}) allows to
prove that the value function of the above problem is the unique
mild solution of the corresponding Hamilton-Jacobi-Bellman
equation (that has the same structure as (\ref{pdeintroduction}).
Moreover the optimal control is expressed in terms of a feedback
that involves the gradient of that same solution to the
Hamilton-Jacobi-Bellman equation. We stress that the usual
application of the Girsanov technique is not allowed (since the
Novikov condition is not guaranteed) and we have to use specific
arguments both to prove the fundamental relation and to solve the
closed loop equation. We adapt some procedure used in \cite{FHT}
to our infinite dimensional framework on infinite horizon.

The paper is organized as follows: the next Section is devoted to
notations; in Section 3 we deal with quadratic BSDEs with random
terminal time; in Section 4 we study the forward backward system
on infinite horizon; in Section 5 we show the result about the
solution to PDE. The last Section is devoted to the application to
the control problem.
\section{Notations}
\label{sec-Notation} The norm of an element $x$ of a Banach space
$E$ will be denoted $|x|_E$ or simply $|x|$, if no confusion is
possible. If $F$ is another Banach space, $L(E,F)$ denotes the
space of bounded linear operators from $E$ to $F$, endowed with
the usual operator norm.

The letters $\Xi$, $H$, $U$ will always denote Hilbert spaces.
Scalar product is denoted $\<\cdot,\cdot\>$, with a subscript to
specify the space, if necessary. All Hilbert spaces are assumed to
be real and separable. $L_2(\Xi,U)$ is the space of
Hilbert-Schmidt operators from $\Xi$ to $U$, endowed with the
Hilbert-Schmidt norm, that makes it a separable Hilbert space. We
observe that if $U=\R$ the space $L_2(\Xi, \rset)$ is the space
$L(\Xi, \rset)$ of bounded linear operators from $\Xi$ to $\R$. By
the Riesz isometry the dual space ${\Xi}^*=L(\Xi, \rset)$ can be
identified with $\Xi$.

By a cylindrical Wiener process with values in
  a Hilbert space $\Xi$, defined on a probability
  space $(\Omega, \calf,\P)$, we mean a family
$\{W_t,\,t\geq 0\}$ of linear mappings from $\Xi$ to
$L^2(\Omega)$, denoted $\xi\mapsto \<\xi, W_t\>$, such that
\begin{description}
  \item[(i)] for every $\xi\in \Xi$,
  $\{ \<\xi,W_t\>,\; t\geq 0\}$ is a real
  (continuous) Wiener process;
  \item[(ii)] for every $\xi_1,\xi_2\in \Xi$ and $t\geq 0$,
  $\E\; (\<\xi_1,W_t\>\cdot \<\xi_2,W_t\>)=
  \<\xi_1,\xi_2\>_\Xi\; t$.
\end{description}

 $(\calf_t)_{t\geq 0}$ will denote, the natural
  filtration of $W$, augmented with the family of
  $\P$-null sets.
The filtration $(\calf_t)$ satisfies the usual conditions. All the
concepts of measurably for stochastic processes refer to this
filtration. By $\calb(\Lambda)$ we mean the Borel $\sigma$-algebra
of any topological space $\Lambda$.



We also recall notations and basic facts on a class of
differentiable maps acting among Banach spaces, particularly
suitable for our purposes (we refer the reader to \cite{FT1} for
details and properties). We notice that the use of G\^ateaux
differentiability in place of Fr\'echet differentiability is
 particularly suitable when dealing with evaluation (Nemitskii)
  type mappings on spaces of summable functions.

Let now $X$, $Z$, $V$ denote Banach spaces. We say that a mapping
$F:X\to V$ belongs to the class $\calg^1 (X,V)$ if it is
continuous, G\^ateaux differentiable on $X$, and its G\^ateaux
derivative $\nabla F:X\to L(X,V)$ is strongly continuous.

 The last requirement is equivalent to the
 fact that for every
$h\in X$ the map $\nabla F(\cdot)h:X\to V$ is continuous. Note
that $\nabla F:X\to L(X,V)$ is not continuous in general if
$L(X,V)$ is endowed with the norm operator topology; clearly, if
this happens then $F$ is Fr\'echet differentiable on $X$. It can
be proved that if $F\in\calg^1(X,V)$ then $(x, h)\mapsto \nabla
F(x)h$ is
  continuous from $X\times X$ to $V$; if, in addition,
$G$ is in $\calg^1(V,Z)$ then $G( F)$ belongs to $\calg^{1} (X,Z)$
and
   the chain rule holds:
   $\nabla(G( F))(x)=\nabla G(F(x))\nabla F(x)$.

When $F$ depends on additional arguments,
 the previous definitions and properties  have
obvious generalizations.

\section{Quadratic BSDEs with random terminal time}


Let $\tau$ be an $\m F_t$-stopping time. It is not necessarily
bounded or $\P$-a.s. finite. We work with a function $F$ defined
on $\Omega \times [0, \infty) \times \R \times \Xi^*$ which takes
its values in $\R$ and such that $F(\cdot,y,z)$ is a progressively
measurable process for each $(y,z)$ in $\rset \times \Xi^*$. We
define the following sets of $\m F_t$-progressively measurable
processes $(\psi_t)_{t \geq 0}$ with values in a Hilbert space
$K$:
$$\m M^{2, -2 \lambda}(0,\tau;K)=\left\{ \psi: \E\left(\int_0^{\tau} e^{-2 \lambda s}|\psi_s|^2ds\right)< \infty \right\}, $$
 $$\m M_{loc}^2(0,\tau;K) =\left\{ \psi: \E\left(\int_0^{t \wedge \tau} |\psi_s|^2ds\right)< \infty \quad \forall t \geq 0\right\}. $$
We want to construct an adapted process $(Y,Z)_{t \geq 0}$ which
solves the BSDE
\begin{equation}\label{bsde-stime}
   -dY_t= \textbf{1}_{t \leq \tau} (F(t,Y_t,Z_t)dt-  Z_t dW_t), \quad Y_{\tau}=\xi \mbox{ on } \{\tau<\infty\}.
\end{equation}

We assume that:
\begin{hyp}
\label{hdiff} There exist $C\geq 0$ and $\alpha\in(0,1)$ such that
\begin{enumerate}
\item $|F(t,y,z)| \leq C\left( 1+|y| + |z|^2 \right)$; \item
$F(t,\cdot,\cdot)$ is $\mathcal{G}^{1,1}(\R\times L_2(\Xi,
\R);\R)$;
\item $\left| \nabla_z F(t,y,z) \right| \leq C\left( 1  + |z|
\right)$; \item $\left| \nabla_y F(t,y,z) \right| \leq C\left( 1 +
|z| \right)^{2\alpha}$.

Moreover we suppose that there exist two constants $K\geq 0$ and
$\lambda>0$ such that $d\P\otimes dt$ a.e.: \item $F$ is monotone
in $y$ in the following sense:
$$\forall y,y' \in \rset, z \in \Xi^*, \quad <y-y',F(t,y,z)-F(t,y',z)> \leq -\lambda|y-y'|^2;$$
\item $|F(t,0,0)| \leq K$; \item $\xi$ is a $\m
F_{\tau}$-measurable bounded random variable; we denote by $M$
some real such that $|\xi| \leq M$ $\P$-a.s.

\end{enumerate}
\end{hyp}

We call solution of the equation a pair of progressively
measurable processes $(Y_t,Z_t)_{t \geq 0}$ with values in $\R
\times \Xi^*$ such that
\begin{enumerate}
   \item $Y$ is a bounded process and $Z \in \m M_{loc}^{2}(0,\tau; \Xi^*)$;
\item On the set $\{\tau <\infty \}$, we have $Y_{\tau}=\xi$ and
$Z_t=0$ for $t> \tau$; \item $\forall T\geq 0$, $\forall t \in
[0,T]$ we have $Y_{t\wedge \tau}= Y_{T\wedge \tau} +\int_{t\wedge
\tau}^{T\wedge \tau} F(s,Y_s,Z_s)ds - \int_{t\wedge \tau}^{T\wedge
\tau} Z_s dW_s$.
\end{enumerate}

Before giving the main result of this section we prove a lemma
which we use in the sequel. The proof involves the Girsanov
transform and results of the bounded mean oscillation (BMO, for
short) martingales theory.

Here we recall a few well-known facts from this theory following
the exposition in \cite{Kaz94}.
Let M be a continuous local $(P ,\m F)$-martingale satisfying $M_0
= 0$. Let $1 \leq  p < \infty$. Then M is in the normed linear
space BMOp if
      $$||M||_{BMO_p}= \sup_{\tau} \left|\left|\E[|M_{T} - M_{\tau}|^p| \m F_{\tau}]^{1/p}\right|\right|_{\infty}<\infty,$$
where the supremum is taken over all stopping time $\tau \leq T$.
By Corollary 2.1 in \cite{Kaz94}, M is a BMOp-martingale if and
only if it is a BMOq-martingale for every $q \geq 1$. Therefore,
it is simply called a BMO-martingale. In particular, M is a
BMO-martingale if and only if
             $$||M||_{BMO_2}= \sup_{\tau} \left|\left|\E[\langle M\rangle_{T} - \langle M \rangle_{\tau}|\m F_{\tau}]^{1/2}\right|\right|_{\infty}<\infty,$$
where the supremum is taken over all stopping time $\tau \leq T$;
$\langle M \rangle$ denotes the quadratic variation of $M$. This
means that local martingales of the form $M_t =\int_0^t \xi_s
dW_s$ are BMO-martingales if and only if
 $$||M||_{BMO_2}= \sup_{\tau} \left|\left|\E \left[\int_{\tau}^T||\xi_s||^2ds \Big|\m F_{\tau}\right]^{1/2}\right|\right|_{\infty}<\infty.$$

The very important feature of BMO-martingales is the following
(see Theorem 2.3 in \cite{Kaz94}): the exponential martingale
$$
\m E(M)_t=\m E_t = \exp\left(M_t - \frac{1}{2} \langle M \rangle_t
\right) \quad  0\leq t \geq T
$$
is a uniformly integrable martingale.

\begin{lemme}\label{LemmaGir}
    Let $(U,V)$, be solutions to
    \begin{equation}\label{eqlemma}
        U_t = \xi + \int_t^T \ind_{s\leq \tau} [a_sU_s +b_sV_s+\psi_s\,] ds - \int_t^T V_s\, dW_s
    \end{equation}
    where $\xi$ is $\m F_\tau$--measurable and bounded and $a_s,b_s,\psi_s$ are processes such that
\begin{itemize}
   \item[1)] $a_s \leq - \lambda$ for some $\lambda >0$;
\item[2)] $\int_0^{\cdot} b_s dW_s$ is a BMO-martingale; \item[3)]
$|\psi_s| \leq \rho(s)$ where $\rho$ is a deterministic function.
\end{itemize}
Moreover we assume that $U$ is bounded. Then we have $\P$-a.s. for
all $t \in [0,T]$
    \begin{equation*}
        |U_t| \leq e^{-\lambda(T-t)}\| \xi\|_\infty  + \int_t^T \rho(s) e^{-\lambda (s-t)}\, ds.
    \end{equation*}
\end{lemme}

\begin{proof}

Let $(U,V)$ be a solution of the BSDE  (\ref{eqlemma}) such that
$U$ is bounded.

We fix $t \in \rset_+$ and set for $s \geq t$ $e_s = e^{\int_{t\wedge\tau
}^{s\wedge\tau} a_r\, dr}$. By Ito's formula
we have,
$$
  U_{t} = e_T\xi + \int_t^T \textbf{1}_{s\leq \tau}e_s \psi_sds- \int_t^T e_s V_s(dW_s - b_s).
$$
Let $\Q_T$ the probability measure on $(\Omega, \m F_T)$ whose
density with respect to $\P_{| \m F_T}$ is
$$\m E_T  =\exp \left(\int_0^T b_sdW_s - \frac{1}{2} \int_0^T |b_s|^2ds \right).$$
By assumption $\int_0^{\cdot} b_s dW_s$ is a BMO-martingale and
the probability measures $\Q_T$ and $\P_{| \m F_T}$ are mutually
absolutely continuous and $\bW_t = W_t - \int_0^t b_r\, dr$ for $0
\leq t \leq T$ is a Brownian motion under $\Q_T$.

Taking the conditional expectation with respect to $\m F_t $ we
get
$$| U_t| \leq  \E^{\Q_T}\left[ e_T |\xi|+ \int_t^T e_s|\psi_s|ds \: \Big| \calf_t\right], \quad \Q_T \mbox{ a.s.}$$
and thanks to 3)
$$
|U_t| \leq (\m E_t)^{-1} \e \left( \m E_T e_T |\xi|+ \int_t^T
\rho(s)e_s  ds\: \Big|\: \m F_{t} \right).
$$
But from 1) $a_s\leq -\lambda$ and, for all $s \geq t$ $ e_s \leq
e^{-\lambda(s- t)} $ $\P$-a.s., from which we get $\P$-a.s.
$\forall t \in[0,T]$
\begin{equation*}
|U_t| \leq e^{-\lambda(T-t)}||\xi||_{\infty} +\int_t^T  \rho(s)
e^{-\lambda(s-t)} ds.
\end{equation*}
\end{proof}

\begin{cor}\label{corGir}
    Let $(Y^i,Z^i)$, $i=1,2$, be solutions to
    \begin{equation*}
        Y^i_t = \xi^i + \int_t^T \ind_{s\leq \tau} F^i(s,Y^i_s,Z^i_s)\, ds - \int_t^T Z^i_s\, dW_s
    \end{equation*}
    where $\xi^i$ is $\m F_\tau$--measurable and bounded. We assume that $Y^1$ and $Y^2$ are bounded and that the $Z^i$ are such that
$\int_0^{\cdot} Z^i_d dW_t$ are BMO-martingales. Moreover $F^1$ is
$-\lambda$-monotone in the following sense: there exists $\lambda
>0$ such that
\begin{equation*}\label{F1monotone}\forall y,y' \in \rset, z \in \Xi^*, \quad
<y-y',F^1(t,y,z)-F^1(t,y',z)> \leq
-\lambda|y-y'|^2;\end{equation*} and verifies
    \begin{equation*}\label{F1}
        |F^1(t,y,z)-F^1(t,y,z')| \leq C\,|z-z'| \left(1+|z|+|z'|\right).
    \end{equation*}
    We assume moreover that
    \begin{equation*}\label{F1-F2}
        |F^1(t,Y^2_t,Z^2_t)-F^2(t,Y^2_t,Z^2_t)| \leq \rho(t)
    \end{equation*}
    where $\rho$ is a deterministic function. Then we have
    $\P$-a.s. for all $t \in [0,T]$
    \begin{equation*}
        |Y^1_t-Y^2_t| \leq e^{- \lambda (T-t)}\| \xi^1-\xi^2 \|_\infty  + \int_t^T \rho(s) e^{-\lambda (s-t)}\,
        ds.
    \end{equation*}
\end{cor}

\begin{proof}

Let $(Y^1,Z^1)$ and $(Y^2,Z^2)$ be solutions of the BSDE with data
respectively $(\xi^1,F^1)$ and $(\xi^2,F^2)$ such that $Y^1$ and
$Y^2$ are bounded. We set $\by= Y^1-Y^2$ and $\bz=Z^1-Z^2$. It is
enough to write the equation for the difference $\by=Y^1 - Y^2$
$$d\by_t= -\textbf{1}_{t\leq \tau}[F^1(t,Y_t^1,Z_t^1)- F^2(t,Y_t^2,Z_t^2)dt + \bz_t dW_t]$$
as
$$d\by_t= -\textbf{1}_{t\leq \tau}[(a_t \by_t + b_t \bz_t + \psi_t)dt + \bz_t dW_t].$$
using a linearization procedure by setting
\begin{equation*}
a_s=  \begin{cases}\dis{\frac{F^1(s,Y_s^1,Z_s^1)-F^1(s,Y_s^2,Z_s^1)}{Y_s^1-Y_s^2}}, \mbox{ if } Y_s^1-Y^2_s\neq 0\\
          -\lambda \mbox{ \qquad \qquad \qquad \qquad \qquad \qquad    otherwise }
           \end{cases}
\end{equation*}

\begin{equation*} b_s=\begin{cases}\dis{ \frac{F^1(s,Y_s^2,Z_s^1)-F^1(s,Y_s^2,Z_s^2)}{|Z_s^1-Z_s^2|^2} (Z_s^1-Z^2_s)},  \mbox{ if } Z_s^1-Z^2_s\neq 0\\
0 \mbox{ \qquad \qquad \qquad \qquad \qquad \qquad \qquad \qquad
\qquad      otherwise .}
          \end{cases}
\end{equation*}
and
\begin{equation*}
\psi_s= F^1(s,Y_s^2,Z^2_s) -F^2(s,Y_s^2,Z^2_s)\qquad\qquad
\end{equation*}

\end{proof}

Now we can state the main result of this section, concerning the
existence and uniqueness of solutions of BSDE (\ref{bsde-stime}).
\begin{thm}\label{E-UquadBSDEIH}
Under assumption \ref{hdiff} there exists a unique solution
$(Y,Z)$ to BSDE (\ref{bsde-stime}) such that $Y$ is a continuous
and bounded process and $Z$ belongs to $\m
M_{loc}^{2}(0,\tau;\Xi^*)$.
\end{thm}

\begin{proof}{\bf Existence.} We adopt the same strategy as in \cite{BH98}
and \cite{Ro}, with some significant modifications.

Denote by $(Y^n,Z^n)$ the unique solution to the BSDE
\begin{equation}\label{eq-n1}
   Y_t^n= \xi\textbf{1}_{\tau \leq n}+\int_t^n \textbf{1}_{s \leq \tau}F(s,Y_s^n,Z_s^n) ds - \int_t^n Z_s^n dW_s, \qquad 0 \leq t\leq n.
\end{equation}

We know from results of \cite{Kob00} that under
\ref{hdiff}-1,2,3,4 the BSDE (\ref{eq-n1}) has a unique bounded
solution and that
\begin{equation*}
\label{equbmo1} \left\| \sup\nl_{t \in[0,\tau\wedge n]} |Y_t|^n
\right\|_\infty \leq (||\xi||_{\infty}+Cn) e^{Cn}\end{equation*}
and there exists a constant $C=C_n$, which depends on $\left\|
\sup\nl_{t\in[0,\tau\wedge n]} \left| Y_t^n \right|
\right\|_\infty$, such that
\begin{equation*}
\label{equbmo2} \left\| \int_0^\cdot Z_s^n\cdot
dW_s\right\|_{BMO_2} \leq C_n.
\end{equation*}
Now we study the convergence of the sequence of processes
$(Y^n,Z^n)$.

(i) First of all we prove that, thanks to the assumptions of
boundedness and monotonicity \ref{hdiff}-5,6, $Y^n$ is a process
bounded by a constant independent on $n$. Applying the Corollary
\ref{corGir} we have that
 $\P$-a.s. $\forall n \in \N$, $\forall t\in [0,n]$
\begin{equation}\label{Y_nbounded}
 |Y_{t}^n| \leq  e^{- \lambda(n-t)}||\xi\textbf{1}_{\tau \leq n}||_{\infty} +\int_t^n e^{- \lambda(s-t)}|F(s,0,0)|\,ds \leq M+ \frac{K}{\lambda}.
\end{equation}

Moreover we can show that for each $\epsilon >0$
\begin{equation}\label{supZn}\sup_{n \geq 1}\E (\int_0^{\tau} e^{-  \epsilon
s}|Z^n_s|^2ds)< \infty .\end{equation} To obtain this estimate we
take the function $\varphi(x)= \left(e^{2Cx}-2Cx-1\right)/(2C^2)$
which has the following properties:
$$\varphi'(x) \geq 0 \mbox{ if } x\geq 0,$$
$$\frac{1}{2} \varphi''(x) - C \varphi'(x)=1 .$$
Thanks to (\ref{Y_nbounded}) we can say that there exist a
constant $K_0$ such that $ \forall s \in [0,T]$, $Y_s^n+K_0 \geq
0$, $\P$-a.s. Now, if we calculate the Ito differential of $e^{-
\epsilon t}\varphi(Y^n_t +K_0)$, using the previous properties, we
have (\ref{supZn}).

(ii) Now we prove that the sequence $(Y^n_t)_{n\geq 0}$ converges
almost surely. We are going to show that it is an almost definite
Cauchy sequence.

We define $Y^n$ and $Z^n$ on the whole time axis by setting
$$Y^n_t= \xi\textbf{1}_{\tau \leq n}, \quad Z^n_t=0, \quad \mbox{ if }  t>n.$$

Fix $t\leq n\leq m$ and set $\wy=Y^m-Y^n$, $\wz=Z^m-Z^n$ and
$\widehat{F}(s,y,z)=\mathbf{1}_{s \leq n}F(s,y,z)$. We get, from
Ito's formula
$$\wy_t= \wy_ m+\int_{t}^{m} \textbf{1}_{s\leq \tau}(F(s,Y_s^m,Z_s^m)- \widehat{F}(s,Y_s^n,Z^n_s))ds - \int_{t}^{m} \wz_s dW_s.$$
We note that
$$|F(s,Y_s^n,Z_s^n)- \widehat{F}(s,Y_s^n,Z^n_s))|=|\mathbf{1}_{s>n}F(s,\xi \textbf{1}_{\tau \leq n},0)| \leq C(1+ M)\mathbf{1}_{s>n}.$$
Hence, we can apply the Corollary \ref{corGir} with $\xi^1=\xi
\textbf{1}_{\tau \leq m}$ and $\xi^2=\xi \textbf{1}_{\tau \leq
n}$, $F^1=F$ and $F^2=\widehat{F}$,
$\rho(t)=C(1+M)\mathbf{1}_{s>n}$ and state that $\forall n,m \in
\nset$, with $n \leq m$ and $\forall t \in [0,n]$, $\P$-a.s.
\begin{multline}\label{Ycauchy}
|Y_{t}^m -Y_{t}^n| \leq e^{-\lambda(m-t)}||\xi \textbf{1}_{\tau
\leq m}- \xi \textbf{1}_{\tau \leq n}||_{\infty}+\int_{n}^{m}
C(1+M)e^{-\lambda (s-t)}\,ds \leq \\
\leq \left(M+\frac{C(1+M)}{\lambda}\right)e^{- \lambda (n- t)}.
\end{multline}
The previous inequality implies that for each $t \geq 0$ the
sequence of random variable $Y_t^n$ is a Cauchy sequence in
$L^{\infty}(\Omega)$, hence converges to a limit, which we denote
$Y_t$. If $m$ goes to infinity in the last inequality, it comes
that $\P$-a.s., $\forall \, 0\leq t \leq n$
\begin{equation}\label{y-yn}|Y_{t}^n-Y_{t}| \leq \beta e^{-\lambda(n-t)}, \qquad\mbox{
where   } \beta=M+\frac{C(1+M)}{\lambda}.\end{equation} This
inequality implies that the sequence of continuous processes
$(Y^n)_{n \in \nset}$ converges almost surely to $Y$ uniformly
with respect to $t$ on compact sets. The limit process $Y$ is also
continuous and from (\ref{Y_nbounded}) we have that $\forall t \in
\rset_{+}$ $|Y_{t}| \leq M + \frac{K}{\lambda}$.

(iii) We show that the sequence $(Y_n)_n$ also converges in the
space $\m M^{2,- 2 \lambda}(0, \tau; \R)$. Indeed we have
$$\E \left[\int_0^{\tau} e^{-2\lambda t} |Y_t^n -Y_t|^2dt\right]= \E \left[\int_0^{n\wedge \tau} e^{-2\lambda t} |Y_t^n -Y_t|^2dt\right]
+\E\left[\int_{n\wedge \tau}^{\tau} e^{-2\lambda t} |Y_t^n
-Y_t|^2dt\right]$$ and using the inequality (\ref{y-yn}) for the
first term, we get that
$$\E \left[\int_0^{n \wedge \tau} e^{-2\lambda t} |Y_t^n -Y_t|^2dt\right]\leq   {\beta}^2n e^{-2 \lambda
n}.$$ In addition, from the definition of $Y^n_t$ on $\rset_+$, we
know that $\forall t > n$ $Y^n_t=\xi \textbf{1}_{\tau \leq n}$.
Hence, \begin{multline*}\E \left[\int_{n \wedge \tau}^{\tau}
e^{-2\lambda t} |Y_t^n -Y_t|^2dt\right]= \E\left[\textbf{1}_{n <
\tau}\int_{n}^{ \tau}e^{-2 \lambda t}|Y_t-\xi\textbf{1}_{n <
\tau}|^2 dt\right]\leq \\\leq 4 \E \left[\textbf{1}_{n
<\tau}\left(M+ \frac{K}{\lambda}\right)^2\int_{n}^{ \tau}e^{-2
\lambda t}dt \right] \leq \frac{2}{\lambda}\left(M+
\frac{K}{\lambda}\right)^2e^{-2 \lambda n}.\end{multline*} Finally
we have
$$\E \left[\int_0^{\tau} e^{-2\lambda t} |Y_t^n -Y_t|^2dt\right]\leq  e^{-2 \lambda n}\left(n{\beta}^2+
\frac{2}{ \lambda} \left(M+ \frac{K}{\lambda}\right)^2\right).$$
Hence $(Y^n)$ converges to $Y$ in $\m M^{2,- 2 \lambda}(0, \tau;
\R)$.

(iv) To continue, we show that the sequence $(Z_n)_n$ is a Cauchy
sequence in the space $ \m M^{2,- 2(\lambda+ \epsilon)}(0,
\tau;\Xi^*)$.

Fix $t\leq n\leq m$ and set, as before, $\wy=Y^m-Y^n$,
$\wz=Z^m-Z^n$ and $\widehat{F}(s,y,z)=\mathbf{1}_{s \leq
n}F(s,y,z)$. We write
\begin{equation*}\label{linearization} F(s,Y_s^m,Z_s^m)- \widehat{F}(s,Y_s^n,Z^n_s)= a_s^{n,m} \wy_s +b_s^{n,m}\wz_s + \mathbf{1}_{s>n}F(s,\xi \textbf{1}_{\tau \leq n},0)
\end{equation*}
where
\begin{equation*}
a_s^{n,m}=  \begin{cases}\dis{\frac{F(s,Y_s^m,Z_s^m)-F(s,Y_s^n,Z_s^m)}{Y_s^m-Y_s^n}}, \mbox{ if } Y_s^m-Y^n_s\neq 0\\
          -\lambda \mbox{ \qquad \qquad \qquad \qquad \qquad \qquad    otherwise }
           \end{cases}
\end{equation*}
\begin{equation*} b_s^{n,m}=\begin{cases}\dis{} \frac{F(s,Y_s^n,Z_s^m)-F(s,Y_s^n,Z_s^n)}{|Z_s^m-Z_s^n|^2}(Z_s^m-Z_s^n),  \mbox{ if } Z_s^m-Z^n_s\neq 0\\
0 \mbox{ \qquad \qquad \qquad \qquad \qquad \qquad    otherwise .}
          \end{cases}
\end{equation*}
From Ito's formula we get
\begin{multline*}
|\wy_0|^2 +\int_0^{\tau\wedge m} e^{-2 (\lambda+\epsilon)
s}|\wz_s|^2 \,ds  +\int_0^{\tau \wedge m} 2e^{-2
(\lambda+\epsilon) s} \wy_s \wz_s dW_s  =\\
= e^{-2 (\lambda+\epsilon)\tau \wedge m}|\wy_{\tau \wedge
m}|^2+\int_0^{\tau \wedge m} e^{-2 (\lambda+\epsilon)
s}2(\lambda+\epsilon) |\wy_s|^2 ds+\\
+\int_0^{\tau \wedge m} 2e^{-2 (\lambda+\epsilon) s}
\wy_s[a_s^{n,m} \wy_s +b_s^{n,m}\wz_s ] ds  + \int_{\tau \wedge
n}^{\tau \wedge m} 2e^{-2 (\lambda+\epsilon) s} \wy_s
F(s,\xi\mathbf{1}_{\tau \leq n},0) ds
\end{multline*}
and taking the expectation we have
\begin{multline*}
 \E \int_0^{\tau\wedge m} e^{-2 (\lambda+ \epsilon) s}|\wz_s|^2 ds  \leq   \E e^{-2 (\lambda+\epsilon)\tau \wedge m}|\wy_{\tau \wedge
m}|^2+ \E\int_0^{\tau \wedge m} e^{-2 (\lambda+\epsilon)
s}2\epsilon |\wy_s|^2 ds+\\
+\E \int_0^{\tau\wedge m} 2 e^{-2 (\lambda+ \epsilon)
s}\wy_sb_s^{n,m}\wz_s ds +
  \E\int_{\tau\wedge n}^{\tau\wedge m} 2e^{-2 (\lambda+\epsilon) s} \wy_s F(s,\xi\mathbf{1}_{\tau \leq n},0) ds .
\end{multline*}
Using the fact that $$2 e^{-2 (\lambda+ \epsilon) s}\wy_
sb_s^{n,m}\wz_s \leq 2 |\wy_ s|^2e^{-2 (\lambda
+\epsilon)s}|b_s^{n,m}|^2 + \frac{1}{2}e^{-2 (\lambda+\epsilon)
s}|\wz_s|^2
$$ we get
\begin{multline*}
\E\int_0^{\tau \wedge m} e^{-2 (\lambda+ \epsilon)s}|\wz_s|^2 \,ds
\leq  2\E e^{-2 (\lambda+\epsilon)\tau \wedge m}|\wy_{\tau \wedge
m}|^2+ 2\E\int_0^{\tau \wedge m} e^{-2 (\lambda+\epsilon)
s}2\epsilon |\wy_s|^2 ds+\\
+4 \E \int_0^{\tau \wedge m} |\wy_ s|^2\,e^{-2 (\lambda+
\epsilon)s} |b_s^{n,m}|^2ds  + \E\int_{\tau \wedge n}^{\tau \wedge
m} 4e^{-2 (\lambda+ \epsilon) s} |\wy_s| |F(s,\xi\mathbf{1}_{\tau
\leq n},0)| ds \leq\\
\leq M^2e^{-2(\lambda+ \epsilon)n}+ {\beta}^2e^{- 2\lambda
n}\left(1+ 4\int_0^{\tau\wedge m}e^{-2 \epsilon s}|b^{n,m}_s|^2ds+
4C(1+M) \E\int_0^{\tau \wedge m} e^{-2 \lambda s} |\wy_s| \right).
\end{multline*}
We note that
$$ |b_s^{n,m}|^2  \leq C(1+  |Z_s^n|^2 + |Z_s^m|^2)$$
and by (\ref{supZn}) $\sup_{n\geq 1} E \int_0^{\tau} e^{-2
\epsilon s}|Z_s^n|^2 ds < \infty$. Finally we obtain
$$\E\int_0^{\tau\wedge m} e^{-2 (\lambda+ \epsilon) s}|\wz_s|^2 \,ds
\leq \beta'(1+n)e^{-2\lambda n}
$$
where $\beta'$ depends on $M,\lambda,K$. Moreover we have that
$$\E \left(\int_{m \wedge \tau}^{\tau} e^{-2 (\lambda + \epsilon) s}|\wz_s|^2 \,ds \right)=0$$
hence
\begin{equation*}
\E \left(\int_0^{\tau} e^{-2 (\lambda + \epsilon)s}|\wz_s|^2 \,ds
\right)\leq \beta'(1+n)e^{-2\lambda n}.
\end{equation*}
Hence $(Z^n)$ is a Cauchy sequence in $\m M^{2, -2 (\lambda+
\epsilon)}(0,\tau;\Xi^*)$ and converges to the process $Z$ in this
space.

(v) It remains to show that the process $(Y,Z)$ satisfies the BSDE
(\ref{bsde-stime}).

We already know that $Y$ is continuous and bounded and $Z$ belongs
to $\m M^{2, -2 (\lambda + \epsilon)}(0,\tau;\Xi^*)$.

By definition $\forall n \in \nset$, $\forall T,t$ such that $0
\leq t \leq T \leq n$ we  have
\begin{equation}\label{passaggio al limite} Y^n_{t\wedge \tau} -Y^n_{T\wedge \tau}= \int_{t\wedge \tau}^{T\wedge \tau} F(s,Y_s^n,Z_s^n) - \int_{t\wedge \tau}^{T\wedge \tau} Z_s^n dW_s.
   \end{equation}

Fix $t$ and  $T$. We shall pass to  the limit in $L^1$ in the
previous equality. The sequence $Y^n_{t\wedge \tau}$ converges
almost surely to $Y_t$ and is bounded by $M+\frac{K}{\lambda}$
uniformly in $n$. From Lebesgue's theorem we get that the sequence
converges to $Y_{t\wedge \tau}$ in $L^1$. Moreover, $
\int_{t\wedge \tau}^{T\wedge \tau} Z_s^n dW_s$ converges in to
$\int_{t\wedge \tau}^{T\wedge \tau} Z_s dW_s$ in $L^2$ since
\begin{equation*}
\E \left(\int_{t \wedge \tau}^{T\wedge \tau} Z_s^n dW_s
-\int_{t\wedge \tau}^{T\wedge \tau} Z_s dW_s \right)^2 \leq e^{2
(\lambda+ \epsilon) T}\E \int_0^{T\wedge \tau} e^{-2 (\lambda+
\epsilon) s}|Z_s^n- Z_s|^2ds.
\end{equation*}
We can note that $\int_{t\wedge \tau}^{T\wedge \tau}
F(s,Y_s^n,Z_s^n)ds$ converges to $\int_{t\wedge \tau}^{T\wedge
\tau} F(s,Y_s,Z_s)ds$ in $L^1$. Indeed
\begin{equation*}
\E \left|\int_{t\wedge \tau}^{T\wedge \tau} F(s,Y_s^n,Z_s^n)ds
-\int_{t\wedge \tau}^{T\wedge \tau} F(s,Y_s,Z_s)ds \right| \leq \E
\int_0^T |F(s,Y_s^n,Z_s^n)ds - F(s,Y_s,Z_s)|ds
\end{equation*}
and, by the growth assumption on $F$, the map $(Y,Z) \rightarrow
F(\cdot,Y,Z)$ is continuous from the space
$L^1(\Omega;L^1([0,T];\rset)) \times L^2(\Omega;L^2([0,T];\Xi^*)$
to $L^1(\Omega;L^1([0,T];\rset))$. (By classical result on
continuity of evaluation operators, see e.g. \cite{AP}). Hence,
passing to the limit in the equation (\ref{passaggio al limite}),
we obtain $\forall t,T$ such that $t \leq T$
\begin{equation*}
   Y_{t\wedge \tau} -Y_{T\wedge \tau}= \int_{t\wedge \tau}^{T\wedge \tau} F(s,Y_s,Z_s) - \int_{t\wedge \tau}^{T\wedge \tau} Z_s dW_s.
\end{equation*}
So to conclude the proof, it only remains to check the terminal
condition. Let $\omega \in \{\tau <\infty\}$, and $n \in \nset$
such that $n \geq \tau(\omega)$. Then
$$
\begin{array}{lll}
\dis{|Y_{\tau}- \xi\textbf{1}_{t \leq 2n}|(\omega)}& = &
\dis{|Y_{n \wedge \tau}- \xi\textbf{1}_{t \leq 2n}|(\omega)\leq
|Y_{n \wedge \tau}-Y_{n \wedge \tau}^{2n}|(\omega) +|Y_{n \wedge \tau}^{2n}- \xi\textbf{1}_{t \leq 2n}|(\omega)\leq}\\
& & \dis{\leq \beta e^{\lambda(n \wedge \tau)}(\omega)e^{-2
\lambda n}+|Y_{n \wedge \tau}^{2n}- \xi\textbf{1}_{t \leq
2n}|(\omega)
\leq \beta e^{-\lambda n }}\\
\end{array}
$$
since $Y^{2n}_{n \wedge \tau}=
Y^{2n}_\tau=Y^{2n}_{2n}=\xi\textbf{1}_{t \leq 2n}$ Then,
$Y_{\tau}= \xi$ $\P$-a.s. on the set $\{\tau < \infty\}$, and the
process $(Y,Z)$ is solution for BSDE (\ref{bsde-stime}).

{\bf Uniqueness.}

Suppose that $(Y^1,Z^1)$ and $(Y^2,Z^2)$ are both solutions of the
BSDE (\ref{bsde-stime}) such that $Y^1$ and $Y^2$ are continuous
and bounded and $Z^1$ and $Z^2$ belong to $\m M^2_{loc}(0, \tau;
\Xi^*)$. It follows directly from the Corollary \ref{corGir} that
$\forall t \geq0$
$$Y_t^1 -Y^2_t =0 \mbox{\qquad \qquad $\P$-a.s.}$$
and then, by continuity, $Y^1=Y^2$.

Applying Ito's formula we have that $d\P\otimes dt$-a.e.
$Z^1_t=Z^2_t$.

\end{proof}

\section{The forward-backward system on infinite horizon}

In this Section we use the previous result to study a
forward-backward system on infinite horizon, when the backward
equation has quadratic generator.

We introduce now some classes of stochastic processes with values
in a Hilbert space $K$ which we use in the sequel.

\begin{itemize}
  \item
$L^p(\Omega;L^2(0,s;K))$ defined for $s \in ]0,+\infty]$ and $p\in
[1,\infty)$, denotes the space of equivalence classes of
progressively measurable processes $\psi:\Omega\times [0,s[
\rightarrow K$, such that
$$|\psi|^p_{L^p(\Omega;L^2(0,s;K))}
=\E\left(\int_0^s |\psi_r|^2_K\; dr\right)^{p/2}.
$$
 Elements of $L^p(\Omega;L^2(0,s;K))$ are
identified up to modification.

\item $L^p(\Omega; C(0,s;K))$, defined for $s \in ]0,+\infty[$ and
$p\in [1,\infty[$,
  denotes the space of progressively measurable processes
$\{\psi_t,\, t\in [0,s]\}$ with continuous paths in $K$, such that
 the norm
$$|\psi|_{L^p(\Omega; C( [0,s];K))}^p
=\E \,\sup_{r\in [0,s]}|\psi_r|_K^p
$$
is finite. Elements of $L^p(\Omega;C( 0,s;K))$ are identified up
to indistinguishability.
\item $L_{\rm
loc}^2(\Omega;L^2(0,\infty;K))$
denotes the space of equivalence classes of progressively
measurable processes $\psi: \Omega \times [0,\infty) \rightarrow
K$ such that
$$\forall t>0 \quad \E \int_0^t |\psi_r|^2dr < \infty.$$
\end{itemize}

Now we consider the It\^o stochastic equation for an unknown
process $\{X_s,s\ge 0\}$
 with values in a Hilbert space $H$:
\begin{equation}\label{fsdeA}
X_s=e^{sA}x+\int_0^s e^{(s-r)A}b(X_{r})dr+\int_0^s
e^{(s-r)A}\sigma dW_{r},\ s\ge 0.
\end{equation}

Our assumptions will be the following:

\begin{hyp}\label{hypothesisSDE1}
(i) The operator $A$ is the generator of a strongly continuous
semigroup $e^{tA}$, $t\ge 0$, in a Hilbert space $H$. We denote by
$m$ and $a$ two constants such that $|e^{tA}|\le me^{at}$ for
$t\ge 0$.

(ii) $b:H\rightarrow H$ satisfies, for some constant $L>0$,
$$|b(x)-b(y)|\le L|x-y|,\ x,y\in H.$$

(iii) $\sigma$ belongs to $L(\Xi,H)$ such that $e^{tA}\sigma \in
L_2(\Xi,H)$ for every $t> 0$, and
\begin{equation*}\label{G1}
|e^{tA}\sigma|_{L_2(\Xi,H)}\le Lt^{-\gamma}e^{at},
\end{equation*}
for some constants $L>0$ and $\gamma\in [0,1/2)$.

(iv) We have $b(\cdot)\in {\cal G}^{1}(H,H)$.

(v) Operators $A+b_x(x)$ are dissipative (that is
    $\langle A y,y\rangle +\langle b_x(x)y,y\rangle \leq 0$ for
    all $x\in H$ and $y\in D(A)$).
\end{hyp}

\begin{rem} We note we need of assumptions $(iv)-(v)$ to obtain a
result of regularity of the process $X$ with respect to initial
condition $x$.
\end{rem}

We start by recalling a well known result on solvability of
equation (\ref{fsdeA})
 on a bounded interval, see e.g. \cite{FT1}.
\begin{prop}\label{lemmadz1}
Under the assumption \ref{hypothesisSDE1}, for every $p\in
[2,\infty)$ and $T>0$ there exists a unique process $X^x\in
L^p(\Omega;C(0,T;H))$ solution of (\ref{fsdeA}). Moreover, for all
fixed $T>0$, the map $x \rightarrow X^x$ is continuous from $H$ to
$L^p(\Omega;C(0,T;H))$.
\begin{equation*}\label{fsdeAestimate}
\mathbb E\sup_{r\in [0,T]}|X_r|^p \le C(1+|x|)^p,
\end{equation*}
for some constant $C$ depending only on $q,\gamma,T,L,a$ and $m$.
\end{prop}

We need to state a regularity result on the process $X$. The proof
of the following lemma can be found in \cite{HT}.
\begin{lemme}\label{lemma-nablaX-bdd}
Under Assumptions \ref{hypothesisSDE1} the map $x\rightarrow X^x$
is G\^ateaux differentiable (that is belongs to
$\mathcal{G}(H,L^p(\Omega,C(0,T;H))$). Moreover denoting by
$\nabla_x X^x $ the partial G\^ateaux derivative, then for every
direction $h\in H$, the directional derivative process $\nabla_x
X^x h,t\in \mathbb R$, solves, $\mathbb{P}-a.s.$, the equation
\begin{equation*}
\nabla_x X^x_{t }h=e^{ t A}h+\int_0^t e^{\sigma A} \nabla_x
F(X^x_{\sigma })\nabla_x X^x_{\sigma }h \, d\sigma,\quad t\in
\mathbb R^+.
\end{equation*}
Finally, $\mathbb{P}$-a.s., $|\nabla_x X^x_{t }h|\le |h|$, for all
$t
>0$.
\end{lemme}

The associated BSDE is:
\begin{equation}\label{bsdex}
Y_t^x=Y_T^x +\int_t^TF(X_\sigma^x,Y_\sigma^x,Z_\sigma^x)
d\sigma-\int_t^T Z_\sigma^x dW_\sigma,\quad 0\le t\le T< \infty.
\end{equation}
Here $X^x$ is the unique mild solution to (\ref{fsdeA}) starting
from $X_0=x$. $Y$ is real valued and $Z$ takes values in $\Xi^*$,
$F:H\times \mathbb R\times \Xi^*\rightarrow \mathbb R$ is a given
measurable function.
\medskip

We assume the following on $F$:

\begin{hyp}
\label{hyFx} There exist $C\geq 0$ and $\alpha\in(0,1)$ such that
\begin{enumerate}
\item $|F(x, y,z)| \leq C\left( 1+|y| + |z|^2 \right)$; \item
$F(\cdot,\cdot,\cdot)$ is $\mathcal{G}^{1,1,1}(H\times \R\times
\Xi^*;\R)$ ;
\item $\left| \nabla_x F(x,y,z) \right| \leq C$; \item $\left|
\nabla_z F(x,y,z) \right| \leq C\left( 1  + |z| \right)$; \item
$\left| \nabla_y F(x,y,z) \right| \leq C\left( 1  + |z|
\right)^{2\alpha}$. \item $\lambda>0$ and $F$ is monotone in $y$
in the following sense:
$$  x \in H, y,y' \ \rset, z \in \Xi^* \quad <y-y',F(x,y,z)-F(x,y',z)> \leq -\lambda|y-y'|^2.$$

\end{enumerate}
\end{hyp}

Applying Theorem \ref{E-UquadBSDEIH}, we obtain:
\begin{prop}\label{E-UquadBSDEIH-x}
Let us suppose that Assumptions \ref{hypothesisSDE1} and
\ref{hyFx} hold.
 Then we have:
 \begin{description}
   \item[$(i)$]For any $x\in H$, there exists a solution $(Y^x,Z^x)$ to the
BSDE (\ref{bsdex}) such that $Y^x$ is a continuous process bounded
by ${K}/{\lambda}$, and $Z\in L_{ \rm
loc}^2(\Omega;L^2(0,\infty;\Xi))$ with $\mathbb{E}\int_0^{\infty}
e^{-2(\lambda+ \epsilon) s}|Z_s|^2ds <\infty$. The solution is
unique in the class of processes $(Y,Z)$ such that $Y$ is
continuous and bounded, and $Z$ belongs to $L_{\rm
loc}^2(\Omega;L^2(0,\infty;\Xi))$.
   \item[$(ii)$] For all $T>0$ and $p\geq 1$,
   the map $x \rightarrow (Y^x\big\vert_{[0,T]},
   Z^x\big\vert_{[0,T]}) $
    is continuous from $H$ to the space $L^p(\Omega;C(0,T;\mathbb{R}))\times
    L^p(\Omega;L^2(0,T;\Xi))$.
 \end{description}
\end{prop}

\begin{proof}
Statement (i) is an immediate consequences of Theorem
\ref{E-UquadBSDEIH}.
 Let us prove (ii). Denoting by $(Y^{n,x},Z^{n,x})$ the unique solution of the following BSDE
   (with finite horizon):
\begin{equation}\label{bsdeppx}
Y_t^{n,x}=
\int_t^nF(X_\sigma^x,Y_\sigma^{n,x},Z_\sigma^{n,x})d\sigma-\int_t^n
Z_\sigma^{n,x} dW_\sigma,
\end{equation}
then, from Theorem \ref{E-UquadBSDEIH}again, $|Y_t^{n,x}|\le
\frac{K}{\lambda}$ and the following convergence rate holds:
\begin{equation*}\label{ratex}
|Y_t^{n,x}-Y_t^x|\le \frac{K}{\lambda}\exp\{-\lambda(n-t)\}.
\end{equation*}
Now, if $x'_m\rightarrow x$ as $m\rightarrow +\infty$ then
\begin{eqnarray*}
|Y^{x'_m}_T-Y^x_T|&\le&|Y^{x'_m}_T-Y^{n,x'_m}_T|+|Y^{n,x}_T-Y^{x}_T|+
|Y^{n,x'_m}_T-
Y^{n,x}_T|\\
&\le& 2
\frac{K}{\lambda}\exp\{-\lambda(n-T)\}+|Y^{n,x'_m}_T-Y^{n,x}_T|.\end{eqnarray*}
Moreover for fixed $n$, $Y^{n,x'_m}_T\rightarrow Y^{n,x}_T$ in
$L^p(\Omega,\mathcal{F}_T,\mathbb{P};\mathbb{R})$ for all $p>1$,
by Proposition 4.2 in \cite{BC} Thus $Y^{x'_m}_T\rightarrow
Y^{x}_T$ in $L^p(\Omega,\mathcal{F}_T,\mathbb{P};\mathbb{R})$.

Now we can notice that
$(Y^{x}\big\vert_{[0,T]},Z^{x}\big\vert_{[0,T]})$ is the unique
 solution of the following BSDE (with finite horizon):
$$
Y_t^{x}= Y_T^{x} +
\int_t^TF(X_\sigma^x,Y_\sigma^{x},Z_\sigma^{x})-\int_t^T
Z_\sigma^{x} dW_\sigma,
$$
and the same holds for
$(Y^{x'_m}\big\vert_{[0,T]},Z^{x'_m}\big\vert_{[0,T]})$. By
similar argument as in \cite{BC} we have
\begin{eqnarray*}
\lefteqn{ \E\left[ \sup_{t \in [0,T]}| Y_t^{x}
-Y_t^{x'_m}|^p\right]^{1 \wedge 1/p} +
\E \left[ \left(\int_0^T |Z_{t}^x-Z_{t}^{x'_m}|\right)^{p/2}\right]^{1 \wedge 1/p} } \\
&\leq& C\, \E \left[\left|Y_T^x-Y_T^{x'_m}
\right|^{p+1}\right]^{\frac{1}{p+1}} +
 \E \left[\left(\int_0^T \left| F(s,X^x_s,Y_s,Z_s) - F(s,X^{x'_m}_s,Y_s,Z_s) \right| ds \right)^{p+1} \right]^{\frac{1}{p+1}}
\end{eqnarray*}
and we can conclude that
$(Y^{x'_m}\big\vert_{[0,T]},Z^{x'_m}\big\vert_{[0,T]}) \rightarrow
(Y^{x}\big\vert_{[0,T]},Z^{x}\big\vert_{[0,T]})$ in
$L^p(\Omega;C(0,T;\mathbb{R}))\times L^p(\Omega;L^2(0,T;\Xi))$.
\end{proof}
$ $

We need to study the regularity of $Y^x$. More precisely, we would
like to show that $Y_0^x$ belongs to ${\cal G}^1(H,\mathbb R)$.

We are now in position to prove the main result of this section.

\begin{thm}\label{theoremregularity}
Under Assumption  the map $x\rightarrow Y_0^x$ belongs to ${\cal
G}^1(H,\mathbb R)$. Moreover $|Y_0^x|+|\nabla_x Y_0^x|\leq c$, for
a suitable constant $c$.
\end{thm}

\begin{proof} Fix $n\ge 1$, let us consider the solution
$(Y^{n,x},Z^{n,x})$ of (\ref{bsdeppx}). Then, see \cite{BC},
Proposition 4.2, the map $x\rightarrow
(Y^{n,x}(\cdot),Z^{n,x}(\cdot))$ is G\^ateaux differentiable from
$H$ to $L^p(\Omega,\;C(0,T;\mathbb{R}))\times
L^p(\Omega;L^2(0,T;\Xi^*))$,
 $\forall p\in (1,\infty)$. Denoting
by $(\nabla_x Y^{n,x}h,\nabla_x Z^{n,x}h)$ the partial G\^ateaux
derivatives with respect to $x$ in the direction $h\in H$, the
processes $\{\nabla_x Y^{n,x}_t h,\nabla_x Z^{n,x}_t h,t\in
[0,n]\}$ solves the equation, $\mathbb{P}-a.s.$,
\begin{eqnarray}\label{derivativeyz}
\nabla_x Y^{n,x}_t h &=& \int_t^n\nabla_x
F(X^x_\sigma,Y^{n,x}_\sigma,
Z^{n,x}_\sigma)\nabla_x X^{n,x}_\sigma h\, d\sigma\nonumber\\
& & +\int_t^n \nabla_y F(X^{x}_\sigma ,Y^{n,x}_\sigma
,Z^{n,x}_\sigma)
\nabla_x Y^{n,x}_\sigma h\, d\sigma\\
& & +\int_t^n\nabla_z
F(X^{x}_\sigma,Y^{n,x}_\sigma,Z^{n,x}_\sigma) \nabla_x
Z^{n,x}_\sigma h\, d\sigma -\int_t^n \nabla_x Z^{n,x}_\sigma h \,
dW_\sigma.\nonumber
\end{eqnarray}
We note that we can write the generator of the previous equation
as
$$\phi_{\sigma}^n(u,v) = \psi_{\sigma}^n + a_\sigma^n u + b^n_\sigma v $$
setting $$\psi_{\sigma}^n = \nabla_x F(X^x_\sigma,Y^{n,x}_\sigma,
Z^{n,x}_\sigma)\nabla_x X^{n,x}_\sigma h $$
$$ a_\sigma^n=
 \nabla_y F(X^{x}_\sigma
,Y^{n,x}_\sigma ,Z^{n,x}_\sigma)
 \quad b^n_\sigma=\nabla_z F(X^{x}_\sigma,Y^{n,x}_\sigma,Z^{n,x}_\sigma).$$

By Assumption \ref{hyFx} and Lemma \ref{lemma-nablaX-bdd}, we have
that for all $x,h\in H$ the following holds $\mathbb{P}$-a.s. for
all $n\in \mathbb{N}$ and all $\sigma\in [0,n]$:
$$\begin{array}{c}
   | \psi_{\sigma}^n| =\Big|\nabla_x F(X^{x}_\sigma,Y^{n,x}_\sigma, Z^{n,x}_\sigma)\nabla_x
X^{x}_\sigma h\Big| \leq C |h|,\\
  a_\sigma^n= \nabla_y F(X^{x}_\sigma
,Y^{n,x}_\sigma ,Z^{n,x}_\sigma) \leq - \lambda \leq 0, \qquad
   |b_\sigma^n|= \Big|\nabla_z F(X^{x}_\sigma,Y^{n,x}_\sigma,Z^{n,x}_\sigma)\Big| \leq C(1 + |Z^{n,x}_\sigma|).
  \end{array}
$$


Therefore $\int_0^{\cdot} Z^{n,x}_\sigma dW_{\sigma}$ is a
BMO-martingale. Hence $\int_0^{\cdot} b_s dW_s$ is also a
BMO-martingale and by Lemma \ref{LemmaGir}, we obtain:
\begin{equation*}\label{estimate-Y}
   \sup_{t\in [0,n]} |\nabla_x Y^{n,x}_t|\le C|h|,\quad \mathbb{P}-\hbox{a.s.};
\end{equation*}
and applying It\^{o}'s formula to $e^{-2\lambda t}|\nabla_x
Y^{n,x}_th|^2$ and arguing as in the proof of Theorem
\ref{E-UquadBSDEIH}, points (iii) and (iv), tanks to the
(\ref{supZn}), we get:
\begin{equation*}\label{estimate-YZ}
    \mathbb E\int_0^{\infty} e^{-2\lambda t}(|\nabla_x Y^{n,x}_t h|^2+|\nabla_x
Z^{n,x}_t h|^2)dt\le C_1|h|^2.
\end{equation*}

Fix $x,h\in H$, there exists a subsequence of $\{(\nabla_x Y^{n,x}
h,\nabla_x Z^{n,x} h,\nabla_x Y^{n,x}_0h): n\in \mathbb{N} \}$
which we still denote by itself, such that $(\nabla_x Y^{n,x}
h,\nabla_x Z^{n,x}h)$ converges weakly to $(U^1(x,h),V^1(x,h))$ in
${\cal M}^{2,-2\lambda}(0,\infty; \rset \times \Xi^*)$ and
$\nabla_x Y^{n,x}_0 h$ converges to $\xi(x,h)\in \mathbb{R}$.

Now we write the equation (\ref{derivativeyz}) as follows:
\begin{eqnarray}\label{derivativeyzforward-bis}
\nabla_x Y^{n,x}_t h&=&\nabla_x Y^{n,x}_0 h -\int_0^t
\nabla_x F(X^{x}_\sigma,Y^{n,x}_\sigma,Z^{n,x}_\sigma)\nabla_x X^x_\sigma hd\sigma\nonumber\\
& & -\int_0^t(
\nabla_y F(X^{x}_\sigma,Y^{n,x}_\sigma,Z^{n,x}_\sigma))\nabla_x Y^{n,x}_\sigma hd\sigma\\
& & -\int_0^t\nabla_z
F(X^{x}_\sigma,Y^{n,x}_\sigma,Z^{n,x}_\sigma) \nabla_x
Z^{n,x}_\sigma hd\sigma\nonumber+\int_0^t \nabla_x Z ^{n,x}_\sigma
hdW_\sigma \nonumber
\end{eqnarray}
and define an other process $U^2_{t}(x,h)$ by
\begin{eqnarray} \label{derivativeyzforward}
\displaystyle U^2_{t}(x,h)&=&\displaystyle\xi(x,h)-\int_0^t
\nabla_x F(X^x_{\sigma},Y^x_\sigma,Z^x_\sigma)\nabla_x X^x_\sigma\,hd\sigma \nonumber\\
& &\displaystyle-\int_0^t( \nabla_y F(X^x_{\sigma},Y^x_\sigma,
Z^x_\sigma)U^1_{\sigma}(x,h)d\sigma \\
& &\displaystyle -\int_0^t\nabla_z
F(X^x_{\sigma},Y^x_\sigma,Z^x_\sigma)
V^1_{\sigma}(x,h)d\sigma+\int_0^t V^1_{\sigma}(x,h)dW_\sigma,
\nonumber
\end{eqnarray}
where $(Y^x,Z^x)$ is the unique bounded solution to the backward
equation (\ref{bsdex}), see Proposition \ref{E-UquadBSDEIH-x}.
Passing to the limit in the equation
(\ref{derivativeyzforward-bis}) it is easy to show that $\nabla_x
Y^{n,x}_t h$ converges to $U^2_{t}(x,h)$ weakly in $L^1(\Omega)$
for all $t >0$.

Thus $U^2_t(x,h)=U^1_t(x,h)$, $\P$-a.s. for a.e. $t \in \R^+$ and
$|U^2_t(x,h)| \leq C|h| $.

Now consider the following equation on infinite horizon
\begin{eqnarray}\label{equationUV}
U(t,x,h)&\!=\!&U(0,x,h)-\int_0^t\nabla_x F(X^x_\sigma,Y^x_\sigma,
Z^x_\sigma)\nabla_x X^x_\sigma hd\sigma\nonumber\\
& &-\int_0^t(\nabla_y F(X^x_\sigma,Y^x_\sigma,Z^x_\sigma))
U(t,x,h)d\sigma\\
& &-\int_0^t\nabla_z
F(X^x_\sigma,Y^x_\sigma,Z^x_\sigma)V(\sigma,x,h) d\sigma+\int_0^t
V(\sigma,x,h)dW_\sigma.\nonumber
\end{eqnarray}
We claim that this equation has a solution.

For each $n \in \nset$ consider the finite horizon BSDE (with
final condition equal to zero):
\begin{eqnarray*}
U_n(t,x,h)&\!=\!&\int_t^n \nabla_x F(X^x_\sigma,Y^x_\sigma,
Z^x_\sigma)\nabla_x X^x_\sigma hd\sigma\nonumber\\
& &+\int_t^n(\nabla_y F(X^x_\sigma ,Y^x_\sigma ,Z^x_\sigma ))
U_n(t,x,h)d\sigma\label{Eq_Un}\\
& &+\int_t^n\nabla_z F(X^x_\sigma,Y^x_\sigma,
Z^x_\sigma)V_n(\sigma,x,h)d\sigma-\int_t^n
V_n(\sigma,x,h)dW_\sigma ,\nonumber
\end{eqnarray*}
By the result in \cite{BC} we know that this equation has a unique
solution $(U_n(\cdot,x,h),V_n(\cdot,x,h))\in
 L^p(\Omega;C(0,n;\mathbb{R}))\times
    L^p(\Omega;L^2(0,n;\Xi^*))$.
The generator of this equation can be rewrite as
$$\phi_t(u,v)=  \psi_t + a_t u+b_t v$$
where $ \psi_t= \nabla_x F(X^x_t,Y^x_t, Z^x_t)\nabla_x X^x_t$ and
$|\psi_t| \leq C|h|$, $a_t =\nabla_y F(X^x_\sigma ,Y^x_\sigma
,Z^x_\sigma )\leq - \lambda$, $b_t=\nabla_z
F(X^x_\sigma,Y^x_\sigma, Z^x_\sigma)$ and $|b_t| \leq C(1
+|Z_t^x|)$. On the interval $[0,n]$ the process $\int_0^{\cdot}
Z^x_s dWs$ is a BMO-martingale. Hence, from the Lemma
\ref{LemmaGir} it follows that $\P$-a.s. $\forall n \in \nset$, $
\forall t \in [0,n]$ $|U^n_t| \leq \frac{C}{\lambda}|h|$ and as in
the proof of existence in the Theorem \ref{E-UquadBSDEIH},
 we can conclude that
\begin{enumerate}
\item for each $t\geq 0$ $U^n(t,x,h)$ is a Cauchy sequence in
$L^{\infty}(\Omega)$ which converges to a process $U$ and
$\P$-a.s., $ \forall t \in [0,n]$
\begin{equation*}\label{StimaUn}|U^n(t,x,h) -U(t,x,h)| \leq \frac{C}{\lambda}|h|e^{-\lambda(n-t)};
\end{equation*}
\item $V^n(\cdot,x,h)$ is a Cauchy sequence in
$L^2_{loc}(\Omega;L^2([0,\infty); \Xi^*)$; \item The processes
limit $(U(\cdot,x,h),V(\cdot,x,h)$ satisfy the BSDE
(\ref{equationUV}).
\end{enumerate}
Moreover still from Lemma \ref{LemmaGir} we get that the solution
is unique.

Coming back to equation (\ref{derivativeyzforward}), we have that
$(U^2(x,h),V^1(x,h))$ is solution in $\mathbb{R}^+$ of the
equation (\ref{equationUV}).

In particular we notice that  $U(0,x,h)=\xi(x,h)$ is the limit of
$\nabla_x Y^{n,x}_0 h$ (along the chosen subsequence). The
uniqueness of the solution to (\ref{equationUV}) implies that in
reality $U(0,x,h)=\lim_{n\to\infty} \nabla_x Y^{n,x}_0 h$ along
the original sequence.

Now let $x_m\rightarrow x$.
\begin{eqnarray}\label{Eq_U_Un}
|U(0,x,h)-U(0,x_m,h)|\leq |U(0,x,h)- U^n(0,x,h)| + |U^n(0,x,h) - U^n(0,x_m,h)| + \\
+|U^n(0,x_m,h)-U(0,x_m,h)| \leq\frac{2 C}{\lambda} e^{-\lambda
n}|h|+ |U_n(0,x,h)-U_n(0,x_m,h)|, \nonumber
\end{eqnarray}
where we have used the (\ref{StimaUn}).
 We now notice that
 $\nabla_x F$, $\nabla_y F$, $\nabla_z F$ are, by assumptions,
 continuous and $|\nabla_x F| \leq C$, $|\nabla_y F| \leq C(1 +| Z|)^{2 \alpha}$, $|\nabla_z F| \leq C(1 +| Z|)$ . Moreover
 the following statements on continuous dependence on $x$ hold:

\noindent maps $x\rightarrow X^x$, $x\rightarrow \nabla_x X^x h$
are continuous from $H\rightarrow L_{\cal P}^p(\Omega;C(0,T;H))$
(see \cite{FT1} Proposition 3.3);

\noindent the map $x\rightarrow Y^x\big\vert_{[0,T]}$ is
continuous from $H$ to $ L_{\cal P}^p(\Omega;C(0,T;\mathbb{R}))$
(see Proposition \ref{E-UquadBSDEIH-x} here);

\noindent the map $x\rightarrow Z^x\big\vert_{[0,T]}$ is
continuous from $H$ to $
    L_{\cal P}^p(\Omega;L^2(0,T;\Xi))$ (see Proposition \ref{E-UquadBSDEIH-x} here ).

We can therefore apply to (\ref{Eq_Un}) the continuity result of
\cite{FT1}
 Proposition 4.3 to obtain in particular that $U_n(0,x'_m,h)\rightarrow
U_n(0,x,h)$ for all fixed $n$ as $m\rightarrow \infty$. And by
(\ref{Eq_U_Un}) we can conclude that $U(0,x'_m,h)\rightarrow
U(0,x,h)$ as $m\rightarrow \infty$.

Summarizing $U(0,x,h)=\lim_{n\to\infty}\nabla_x Y^{n,x}_0 h$
exists, moreover it is
 clearly linear in $h$ and verifies $|U(0,x,h)|\leq C|h|$, finally it is
continuous in $x$ for every $h$ fixed.

Finally, for $t>0$,
$$\begin{array}{rcl}
\displaystyle \lim_{t\searrow 0}\frac{1}{t}[Y^{x+th}_0-Y^{x}_0]
&=& \displaystyle \lim_{t\searrow 0} \frac{1}{t}
\lim_{n\rightarrow +\infty}[Y^{n,x+th}_0-Y^{n,x}_0]=
\lim_{t\searrow 0} \lim_{n\rightarrow +\infty}\int_0^1 \nabla_x Y_0^{n,x+\theta th}h d\theta\\
\displaystyle &=& \displaystyle \lim_{t\searrow 0}
\int_0^1U(0,x+\theta th)h d\theta=U(0,x)h
\end{array}$$
and the claim is proved.
\end{proof}

\section{Mild Solution of the elliptic PDE}
\label{sec-Kolm} Now we can proceed as in \cite{FT3}. Let us
consider the forward equation
\begin{equation}\label{for-system}
X_s=e^{sA}x+\int_0^s e^{(s-r)A}b(X_{r})dr+\int_0^s
e^{(s-r)A}\sigma dW_{r},\ s\ge 0.
\end{equation}
Assuming that Assumption \ref{hypothesisSDE1} holds, we define in
the usual way the transition semigroup
 $(P_t)_{t\geq 0}$, associated to
the process $X$:
\begin{equation*}\label{defsemigruppo}
    P_t[\phi](x)= \E\; \phi(X_t^x),
\qquad x\in H,
\end{equation*}
for every bounded measurable function $\phi:H\to \R$. Formally,
the generator $\call$ of $(P_{t})$ is the operator
$$
\call\phi (x)=\frac{1}{2}{\rm Trace }\left( \sigma
\sigma^*\nabla^2\phi(x)\right) + \< Ax+b(x),\nabla \phi(x)\>.
$$
In this section we address solvability of the non linear
stationary Kolmogorov equation:
\begin{equation}\label{kolmogorovnonlineare}
\call v(x) + F (x,v(x),\nabla v(x)\, \sigma)=0,\qquad x\in H,
\end{equation}
when the coefficient $F$ verifies Assumption \ref{hyFx}.
 Note that, for $x\in H$, $\nabla v(x)$
belongs to $H^*$, so that $\nabla v(x)\, \sigma$ is in $\Xi^*$.

\begin{df}\label{defdisoluzionemild}
We say that a function $v: H\to\R$ is a mild solution of the non
linear stationary Kolmogorov equation (\ref{kolmogorovnonlineare})
if the following conditions hold:
\begin{description}
  \item[(i)]
$v\in\calg^{1}(H,\R)$ and $\exists\, C>0$ such that $ |v(x)|\leq
C$, $
  |\nabla_x v(x)h|\leq C\;|h|,
  $ for all $x, h\in H$;
  \item[(ii)] the following equality holds,
  for every $x\in H$ and $T\geq 0$:
\begin{equation}\label{solmild}
  v(x) =e^{-\lambda T}\;P_{T}[v](x)
  +\int_0^T e^{-\lambda t}\;P_{t}\Big[
F \Big(\cdot,v(\cdot),\nabla v(\cdot)\, \sigma \Big)+ \lambda
v(\cdot) \Big](x)\; dt.
\end{equation}
\end{description}
where $\lambda$ is the monotonicity constant in Assumption
\ref{hyFx}.
\end{df}

Together with equation (\ref{for-system}) we also consider the
backward equation
\begin{equation}\label{back-system}
Y_t -Y_T +\int_t^T Z_s dW_s = \int_t^T F(X_s,Y_s,Z_s)ds \qquad
0\leq t \leq T <\infty
\end{equation}
where $F:H \times \rset \times \Xi^* \rightarrow \rset$ is the
same occurring in the nonlinear stationary Kolmogorov equation.
Under the Assumptions \ref{hypothesisSDE1}, \ref{hyFx},
Propositions \ref{lemmadz1}-\ref{E-UquadBSDEIH-x} give a unique
solution $\{X_t^x,Y_t^x,Z_t^x\}$, for $t \geq 0$, of the
forward-backward system (\ref{for-system})-(\ref{back-system}). We
can now state the following

\begin{thm}\label{main}
Assume that Assumption \ref{hypothesisSDE1}, Assumption \ref{hyFx}
and 
 hold then equation
(\ref{kolmogorovnonlineare}) has a unique mild solution given by
the formula
\begin{equation*}
    v(x) =Y_0^x.
\end{equation*}
where $\{X_t^x,Y_t^x,Z_t^x, t \geq 0\}$ is the solution of the
forward-backward system (\ref{for-system})-(\ref{back-system}).
Moreover the following holds:
\begin{equation*}
Y^x_t=v(X^x_t),\quad Z^x_t=\nabla v(X^x_t)\,\sigma.
\end{equation*}
\end{thm}

{\bf Proof.} Let us recall that for $s \geq 0$, $Y_s^{x}$ is
measurable with respect to $\m F_{[0,s]}$ and $\m F_s$; it follows
that $Y_0^{x}$ is deterministic (see also \cite{ElK97}). Moreover,
as a byproduct of Proposition~\ref{theoremregularity}, the
function $v$ defined by the formula $v(x)=Y_0^{x}$ has the
regularity properties stated in
Definition~\ref{defdisoluzionemild}. The proof that the equality
(\ref{solmild}) holds true for $v$ is identical to the proof of
Theorem 6.1 in \cite{FT3}.

\section{Application to optimal control}\label{sec-control}

We wish to apply the above results to perform the synthesis of the
optimal control for a general nonlinear control system on an
infinite time horizon. To be able to use non-smooth feedbacks we
settle the problem in the framework of weak control problems.
Again we follow \cite{FT3} with slight modifications. We report
the argument for reader's convenience.

As above by $H$, $\Xi$ we denote separable real Hilbert spaces and
by $U$ we denote a Banach space.

For fixed $x_{0}\in H$ an {\it admissible control system} (a.c.s)
is given
 by $(\Omega,{\m F},({\m F}_{t})_{ t\geq 0},{\P}, \{{W}_{t}, t\geq 0\}, u)$ where
 \begin{itemize}

\item $(\Omega,{\m F},\mathbb{P})$ is a complete probability space
and $(\mathcal{F}_{t})_{ t\geq 0}$ is a filtration on it
satisfying the usual conditions.

\item $\{W_{t}: t\geq 0\}$ is a $\Xi$-valued cylindrical Wiener
process
    relatively to the filtration
    $({\cal F}_t)_{ t\geq 0}$ and the probability $\P$.

    \item $u:\Omega\times [0,\infty[\rightarrow U$ is a predictable process
    (relatively to $({F}_{t})_{ t\geq 0}$)
    that satisfies the constraint:
 $u_t\in \mathcal{U}$, $\P$-a.s. for a.e.
 $t\geq 0$,
 where $\mathcal{U}$ is a fixed closed subset of $U$.
 \end{itemize}
  To each a.c.s. we associate the mild solution
  $X\in L^{r}_{\calp}(\Omega;C(0,T;H))$
  (for arbitrary $T>0$ and arbitrary $r\geq 1$)
 of the state equation:
 \begin{equation}\label{equazionestato}
\left\{\begin{array}{l}\dis dX_\tau =\left( AX_\tau+b(X_\tau)
+\sigma r(X_\tau,u_\tau)\right)\; d\tau
+\sigma \; dW_\tau,\qquad \tau\geq 0,\\
\dis X_{0}=x\in H,
\end{array}\right.
\end{equation}
and the cost:
\begin{equation}
J(x,u)=\E\, \int_{0}^{+\infty} e^{-\lambda t} g(X_t,u_t)\; dt ,
    \label{funzionalecosto}
\end{equation}
where $g: H\times U\to\R$. Our purpose is to minimize the
functional $J$ over all a.c.s. Notice the occurrence of the
operator $\sigma$ in the control term: this special structure of
the state equation is imposed by our techniques.

We work under the following assumptions.

\begin{hyp}\label{hypcontrol}
\begin{enumerate}
   \item  The process W is a Wiener process in $\Xi$, defined on a complete probability
        space $(\Omega, \m F, \P)$ with respect to a filtration $(\m F_t )$ satisfying the usual conditions.
       \item $A$, $b$ verify Assumption \ref{hypothesisSDE1}.
\item $\sigma$ satisfies Assumption \ref{hypothesisSDE1} (iii)
with $\gamma=0$;
    \item The set $\m U$ is a nonempty closed subset of $U$ .
      \item The functions $ r:H\times U \rightarrow \Xi$, $g : H\times U \rightarrow  \rset$ are Borel measurable and
 for all $ x \in H$, $r( x, \cdot)$ and $g( x, \cdot) $ are continuous functions from
       $ U$ to $\Xi $  and from $U$ to $\rset$, respectively.
      \item  There exists a constant $C \geq 0$ such that for every  $x, x' \in H$ ,
       $ u \in K $ it holds that
                                  $$|r( x, u) -r( x' , u)| \leq C(1 + |u|)|x - x' |,$$
                                               \begin{equation}\label{r-gro} |r( x, u)| \leq C(1 + |u|),
                                               \end{equation}
                                     \begin{equation}\label{2.8}
0 \leq g( x, u) \leq C(1 + |u|^2 ),
                                     \end{equation}
 \item  There exist $R > 0 $ and  $c > 0$ such that for every $ x\in H$ $u \in U$ satisfying $|u| \geq R$,
                                         \begin{equation}\label{2.10}
                                            g( x, u) \geq c|u|^2.
                                         \end{equation}
  \end{enumerate}
\end{hyp}
    We will say that an $(\m F_t )$-adapted stochastic process $\{u_t , t \geq 0\}$ with values
in $U$ is an admissible control if it satisfies
\begin{equation}\label{adco}\E \int_0^{\infty} e^{-\lambda t} |u_t|^2 dt < \infty.
\end{equation}

This square summability requirement is justified by (\ref{2.10}):
a control process which is not square summable would have infinite
cost.

Now we state that for every admissible control the solution to
(\ref{equazionestato}) exists.
\begin{prop}
 Let $u$ be an admissible control. Then there exists a unique,
continuous, $(\m F_t )$-adapted process $X$ satisfying $\E \sup_{t
\in [0,T ]} |X_t |^2 < \infty$, and $\P$-a.s., $t \in [0, T ]$
$$X_t = e^{tA}x +   \int_0^t  e^{(t-s)A} b(X_s ) ds + \int_0^t  e^{(t-s)A}  \sigma dW_s +     \int_0^t  e^{(t-s)A}
\sigma r(X_s , u_s ) ds.$$
\end{prop}
\begin{proof}
The proof is an immediate extension to the infinite dimensional
case of the Proposition 2.3 in \cite{FHT}.
\end{proof}
By the previous Proposition and the arbitrariness of $T$ in its
statement, the solution is defined for every $t\geq 0$. We define
in a classical way the Hamiltonian function relative to the above
problem: for all $x\in H$, $z\in \Xi^*$,
\begin{equation}\begin{array}{l}
\dis F(x,y,z)=\inf\{g( x,u)+zr(x,u): u\in
\mathcal{U}\}- \lambda y\\
\dis \Gamma(x,y,z)=\{u\in U:g( x,u)+zr(x,u) -\lambda y=
F(x,y,z)\}.
\end{array}\label{definhamilton}
\end{equation}
The proof of the following Lemma can be found in \cite{FHT} Lemma
3.1.
\begin{lemme} The map $F$ is a Borel measurable function from $H \times \Xi^*$ to $\rset$.
 There exists a constant $C > 0$ such that
 \begin{equation}
         -C(1 + |z|^2 ) -\lambda y \leq F(x,y, z) \leq g( x, u) + C|z|(1 + |u|) -\lambda y  \quad       \forall u \in \m U.
 \end{equation}

\end{lemme}

We require moreover that
\begin{hyp}\label{hyFcontrollo}
$F$ satisfies assumption \ref{hyFx} 2-3-4.
\end{hyp}

We notice that the cost functional
 is well defined and $J(x,u)<\infty$
for all $x\in H$ and all a.c.s.

By Theorem \ref{main}, the stationary Hamilton-Jacobi-Bellman
 equation relative to the above stated problem, namely:
 \begin{equation}\label{equazioneHJB}
\call v(x) + F( x,v(x),\nabla v(x)\sigma)=0,\qquad x\in H,
\end{equation}
admits a unique mild solution, in the sense of Definition
\ref{defdisoluzionemild}.

\subsubsection{The fundamental relation}
\begin{prop}\label{fund-relation}
Let $v$ be the solution of (\ref{equazioneHJB}). For every
admissible control $u$ and for the corresponding trajectory $X$
starting at $x$ we have
\begin{equation*}\begin{array}{l} J(x,u) = v(x)+\\
\dis{\E\int_{0}^{\infty}e^{-\lambda t}\bigg(-F(X_t,
\nabla v(X_t)\sigma)-\lambda v(X_t) 
 +
\nabla_{x}v(X_t)\sigma r(X_t , u_t ) +g(X_t,u_t)\bigg) \;
dt.}\end{array}\end{equation*}
\end{prop}
\begin{proof}
We introduce the sequence of stopping times
\begin{equation*}\tau_n = \inf \{ t \in [0, T ] :   \int_0^t |u_s |^2 ds \geq n \},
\end{equation*}
with the convention that $\tau_n = T$ if the indicated set is
empty. By (\ref{adco}), for $\P$-almost every $\omega \in
\Omega$, there exists an integer $N(\omega)$ depending on $\omega$
such that
\begin{equation}\label{tn=T} n \geq N(\omega)       \Longrightarrow      \tau_n (\omega) = T.\end{equation}
Let us fix $u_0 \in K$, and for every $n$, let us define
                              $$u_t^n = u_t 1_{t\leq \tau_n} + u_0 1_{t>\tau_n}$$
and consider the equation
\begin{equation}\label{eqstatoXn}\left\{\begin{array}{l}
dX_t^n = b(X_t^n ) dt + \sigma [dW_t + r( X_t^n , u_t^n ) dt],\quad 0\leq t \leq T\\
 X_0^n = x.
      \end{array}\right.
\end{equation}
 Let us define
                           $$    W_t^n := W_t +  \int_0^t   r(X_s^n , u_s^n ) ds \quad 0 \leq t \leq T.$$
From the definition of $\tau_n$ and from (\ref{r-gro}), it follows
that
\begin{equation}\label{4.1}
     \int_0^T  |r( X_s^n , u_s^n )|^2 ds \leq C \int_0^T   (1 + |u_s^n |)^2 ds \leq C\int_0^{\tau_n}  (1 + |u_s |)^2
 ds + C
 \leq C + Cn.
\end{equation}

Therefore defining
$$\rho_n= \exp{\left(\int_0^T -r( X_s^n , u_s^n ) dW_s -  \frac{1}{2} \int_0^T   |r(X_s^n , u_s^n )|^2 ds\right)} $$
the Novikov condition implies that $\E \rho_n = 1$. Setting
$d\P^n_T = \rho_n d\P_{|\m F_T}$, by the Girsanov theorem $W^n$ is
a Wiener process under $\P^n_T$. Relatively to $W^n$ the equation
(\ref{eqstatoXn}) can be written:
\begin{equation}
  \left\{\begin{array}{l}
dX_t^n = b(X_t^n ) dt + \sigma dW_t^n ,\quad 0\leq t \leq T\\
 X_0^n = x.
      \end{array}\right.
\end{equation}


Consider the following finite horizon Markovian forward-backward
system (with respect to probability $\P^n_T$ and to the filtration
generated by $\{W^n_\tau :\tau\in [0,T]\}$).
\begin{equation}\label{backwardforwardtilde}
\left\{\begin{array}{l}\dis X_\tau^n(x) = e^{\tau A}x+\int_0^\tau
e^{(\tau-s)A} b(X_s^n(x))\; ds +\int_0^\tau e^{(\tau-s)A} \sigma
\; dW_s^n, \quad \tau\geq 0,
\\
\dis Y_\tau^n(x) -v(X_T^n(x)) + \int_t^T Z_s^n(x) dW_s^n =
\int_t^T F(X_s^n(x), Y_s^n(x),Z_s^n(x) ) ds,
  \; 0\leq \tau\leq T,
\end{array}\right.
\end{equation}
and let $(X^n(x),Y^n(x),Z^n(x))$ be its unique solution with the
three processes predictable relatively to the filtration generated
by $\{W^n_\tau :\tau\in [0,T]\}$ and:
$\E^n_T\sup_{t\in[0,T]}|X^n_t(x)|^2<+\infty$, $Y^n(x)$ bounded and
continuous, ${\E}_T^n\int_0^T|Z_t^n(x)|^2dt<+\infty$. Moreover,
Theorem \ref{main} and uniqueness of the solution of system
(\ref{backwardforwardtilde}), yields that
\begin{equation}\label{relcontr}
Y^n_t (x)=v( X^n_t(x)), \qquad Z^n_t(x)=\nabla v(
X^n_t(x))G(X^n_t(x)).
\end{equation}
Applying the It\^o formula to
 $e^{- \lambda t }Y^n_t (x)$,
and restoring the original noise $W$ we get
 \begin{equation}\label{4.3}
\begin{array}{l}\dis
e^{- \lambda \tau_n} Y^n_{\tau_n} (x)= e^{- \lambda T} Y^n_T (x)+
\int_{\tau_n}^T \lambda e^{- \lambda t} Y^n_s (x)ds -
\int_{\tau_n}^T e^{- \lambda s }Z^n_s(x) \;d{W}_s
 \\\dis\qquad
 + \int_{\tau_n}^Te^{- \lambda s }\left[F(X^n_s(x),Y^n_s(x),
  Z_s^n(x))-
Z^n_s(x)r(X^n_s,u_s^n) \right] \;ds.
 \end{array}
\end{equation}
We note that for every $p \in [1, \infty)$ we have
\begin{multline}
            \rho_n^{-p }= \exp\left( p \int_0^T r(X_s^n , u_s^n ) dW_s^n -  \frac{p^2}{2}\int_0^T |r(X_s^n , u_s^n )|^2 ds \right)\\
\cdot\exp\left(\frac{p^2-p}{2}\int_0^T |r( X_s^n , u_s^n )|^2 ds
\right).
\end{multline}
By (\ref{4.1}) the second exponential is bounded by a constant
depending on $n$ and $p$, while the first one has
$\P^n$-expectation, equal to 1. So we conclude that $\E^n
\rho_n^{-p} < \infty$. It follows that
\begin{multline*}
\E \left(\int_0^T e^{-2\lambda t} |Z^n_t(x) |^2 dt \right)^{1/2}  \leq \E^n \left(\int_0^T   \rho_n^{-2} |Z_t^n(x) |^2 dt \right)^{1/2}  \leq \\
\leq (\E^n   \rho_n^{-2} )^{1/2} \E^n \left(\int_0^T |Z_t^n (x)|^2
dt \right)^{1/2}< \infty
\end{multline*}
We conclude that the stochastic integral in (\ref{4.3}) has zero
expectation. Using the identification in (\ref{relcontr}) and
taking expectation with respect to $\P$, we obtain
\begin{equation}\label{4.4} \begin{array}{l}
\dis \E e^{- \lambda \tau_n}  Y^n_{\tau_n} =\dis e^{- \lambda T}\E [v(X^n_T(x))] +\E \int_{\tau_n}^T \lambda e^{- \lambda t} Y^n_s (x)ds +\\
\dis\qquad \qquad  \quad +\E\int_{\tau_n}^Te^{- \lambda s
}\left[F(X^n_s(x),Y^n_s(x),
  Z_s^n(x))-
Z^n_s(x)r(X^n_s(x),u_s^n) \right] \;ds\leq \\
\dis \leq e^{- \lambda T}\E [v(X^n_T(x)] +\E \int_{\tau_n}^T
\lambda e^{- \lambda s} Y^n_s (x)ds + \E\,\int_{\tau_n}^T e^{-
\lambda s} g(X^n_s(x),u^n_s)ds.
\end{array}
\end{equation}
Now we let $n \rightarrow \infty$. By Proposition
\ref{E-UquadBSDEIH-x},
                      \begin{equation}\label{2.9}\sup_{t \geq 0} |Y_t^n|=\sup_{t \geq 0} |v(X^n_t)| \leq \frac{K}{\lambda}; \end{equation}
in particular
$$\E \int_{\tau_n}^T \lambda e^{- \lambda s} Y^n_s (x)ds \leq \E \int_{\tau_n}^T \lambda e^{- \lambda s} \frac{K}{\lambda} ds \leq \E K(T - \tau_n)$$
and the right-hand side tends to 0 by (\ref{tn=T}). By the
definition of $u^n$ and (\ref{2.8}),
\begin{multline}\E \int_{\tau_n}^T  g(X^n_s , u^n_s ) ds=\E \int_0^T 1_{s>\tau_n} g( X_s^n , u_0 ) ds \leq \\
   \leq C \E\int_0^T 1_{s>\tau_n} (1 + |u_0 |^2 )ds \leq C\E (T - \tau_n )
\end{multline}
and the right-hand side tends to 0 again by (\ref{tn=T}). Next we
note that, again by (\ref{tn=T}), for $n \geq N (\omega)$ we have
$\tau_n (\omega) = T$ and
                                $ v(X^n_T ) = v(X^n_{\tau_n} ) = v(X_{\tau_n} ) = v(X_T ).$
 We deduce, thanks to (\ref{2.9}), that
$\E v(X^n_T ) \rightarrow \E v(X_T )$, and from (\ref{4.4}) we
conclude that $$\limsup_{n \rightarrow \infty} \E e^{- \lambda
\tau_n} Y^n_{\tau_n} \leq e^{-\lambda T} \E v(X_T ).$$ On the
other hand, for $n \geq N (\omega)$ we have $\tau_n (\omega) = T$
and
                                    $ e^{- \lambda \tau_n}Y_{\tau_n}^n = e^{- \lambda T}Y^n_T = e^{- \lambda T}v(X^n_T) = e^{- \lambda T}v(X_T).$
Since $Y^n$ is bounded, by the Fatou lemma, $\E e^{-\lambda
T}v(X_T ) \leq \liminf_{n \rightarrow \infty} \E e^{- \lambda
\tau_n}Y_{\tau_n}^n$. We have thus proved that
     \begin{equation}\label{4.5}  \lim _{n \rightarrow \infty} \E e^{- \lambda \tau_n}Y_{\tau_n}^n = e^{-\lambda T}\E v(X_T ).
     \end{equation}
Now we return to backward equation in the system
(\ref{backwardforwardtilde}) and write
\begin{multline*}e^{- \lambda \tau_n}Y_{\tau_n}^n =Y^n_0 + \\
+ \int_0^{\tau_n}  -e^{-\lambda t} F(X^n_t,Y^n_t,Z^n_t)dt +
\int_0^{\tau_n} -\lambda e^{-\lambda t}Y_t^n dt+\int_0^{\tau_n}
e^{-\lambda t} Z^n_tdW_t +\int_0^{\tau_n} e^{-\lambda t} Z^n_t
r(X_t^n , u_t^n )dt
\end{multline*}
Arguing as before, we conclude that the stochastic integral has
zero $\P$-expectation. Moreover, we have $Y_0^n = v(x)$, and, for
$t\leq \tau_n$, we also have $u_t^n = u_t$, $X^n_t = X_t$, $Y^n_t
=v(X_t^n)= v(X_t)$ and $Z^n_t = \nabla_x v(X_t) $. Thus, we obtain
 \begin{multline}\E[ e^{- \lambda \tau_n}Y_{\tau_n}^n] = v(x) + \\
+\E \int_0^{\tau_n}e^{- \lambda t}\bigg(-F(X_t ,v(X_t),\nabla_x
v(X_t) \sigma)- \lambda v(X_t)  + \nabla_x v(X_t) \sigma r( X_t ,
u_t )\bigg) dt
\end{multline}
and  \begin{multline}\E \int_0^{\tau_n} e^{- \lambda t} g(X_t ,
u_t ) dt +\E[ e^{- \lambda \tau_n}Y_{\tau_n}^n]  = v(x) + \\ + \E
\int_0^{\tau_n} e^{- \lambda t} \bigg(-F(X_t ,v(X_t),\nabla_x
v(X_t) \sigma)-\lambda v(X_t) + \nabla_x v(X_t) \sigma r( X_t ,
u_t ) + g(X_t , u_t )\bigg) dt .\end{multline} Noting that $-F(
x,y, z)-\lambda y + z r( x, u) + g( x, u) \geq 0$ and recalling
that $g(x, u) \geq 0$ by (\ref{4.5}) and the monotone convergence
theorem, we obtain for $n \rightarrow \infty$,
 \begin{multline}\E\int_0^T e^{-\lambda t}g( X_t , u_t ) dt +  e^{-\lambda T}\E v(X_T )  = v( x)+\\
+ \E \int_0^T e^{- \lambda t} \bigg(-F(X_t , \nabla_x v(X_t)
\sigma ) -\lambda v(X_t) +\nabla_x v(X_t) \sigma r( X_t , u_t ) +
g(X_t , u_t )\bigg) dt.\end{multline} Recalling that $v$ is
bounded,
  letting
$T\rightarrow \infty$,
 we conclude
$$
     \begin{array}{l}
\dis{ J(x,u)= v(x)+}\\
\dis{\E\int_{0}^{\infty}e^{-\lambda t}\left[-F(X_t,v(X_t),
\nabla v(X_t)\sigma)-\lambda v(X_t) 
 +
\nabla_{x}v(X_t)\sigma r(X_t , u_t ) +g(X_t,u_t)\right] dt .}
 \end{array}
$$
 The above equality is known as the {\it fundamental relation} and
 immediately implies that $ v(x)\leq J(x,u)$
  and
 that the equality holds if and only if the following feedback law holds $\P$-a.s. for almost every $t \geq0$:
$$F(X_t,v(X_t),\nabla_x v(X_t)\sigma)=\nabla_x v(X_t)\sigma + g(X_t,u_t) -\lambda v(X_t) $$
where $X$ is the trajectory starting at $x$ and corresponding to
control $u$.
\end{proof}
\subsubsection{Existence of optimal controls: the closed loop equation.}

Next we address the problem of finding a weak solution to the
so-called closed loop equation. We have to require the following

\begin{hyp}\label{existence-minimum}
$\Gamma(x,y,z)$, defined in \ref{definhamilton}, is non empty for
all $x\in H$ and $z\in \Xi^*$.
\end{hyp}

By simple calculation (see \cite{FHT} Lemma 3.1), we can prove
that this infimum is attained in a ball of radius $C(1  + |z|)$,
that is,
$$     F(x,y, z) =            \min_{ u \in \m U,|u| \leq C(1+|z|)} [g(x, u) + z · r(x, u)] -\lambda y,
                         \quad       x \in H , \, y \in \rset,\,    z \in \Xi^* ,$$
and
            \begin{equation}\label{3.4} F( x,y, z) < g(x, u) + z · r(x, u) -\lambda y \quad \mbox{ if } |u| > C(1  + |z|).
            \end{equation}
Moreover, by the Filippov Theorem (see, e.g., [1, Thm. 8.2.10, p.
316]) there exists a measurable selection of $\Gamma$, a Borel
measurable function $\gamma :  H \times\Xi^*\rightarrow \m U$ such
that
\begin{equation}\label{5.2}
 F(x, y,z) = g( x, \gamma( x, z)) + z · r(x, \gamma( x, z)) -\lambda y,      \quad  x \in H ,\, y \in \rset,\,  z \in \Xi^* .\end{equation}
By (\ref{3.4}), we have
\begin{equation}\label{5.3}
|\gamma( x, z)| \leq C(1  + |z|).
      \end{equation}
We define
$$ \underline{u}(x) = \gamma(x, \nabla_xv(X_t)\sigma) \quad \textit{ $\P$-a.s. for a.e } t\geq 0.$$

                     The closed loop equation is
\begin{equation}\label{closedloop}
\left\{\begin{array}{l}
 dX_t = AX_t dt + b(X_t ) dt + \sigma [dWt + r( X_t , \underline{u}(X_t )) dt]\quad t \geq 0\\
             X_0 = x
 \end{array}\right.
\end{equation}
By a weak solution we mean a complete probability space
$(\Omega,\m F, \P)$ with a filtration $(\m F_t )$ satisfying the
usual conditions, a Wiener process $W$ in $\Xi$ with respect to
$\P$ and $(\m F_t )$, and a continuous $(\m F_t )$-adapted process
$X$with values in $H$ satisfying, $\P$-a.s.,
$$\int_0^{\infty}e^{-\lambda t}|\underline{u}(X_t)|^2dt < \infty$$
and such that (\ref{closedloop}) holds. We note that by
(\ref{r-gro}) it also follows that
                            $$ \int_0^{} |r( X_t , \underline{u}( X_t ))|^ 2dt < \infty, \quad    \P-a.s.,$$
so that (\ref{closedloop}) makes sense.

\begin{prop}\label{solutionclosedloop} Assume that $b, \sigma, g$ satisfy Assumption \ref{hypcontrol}, $F$ verifies Assumption \ref{hyFcontrollo} and Assumption \ref{existence-minimum} holds. Then
 there exists a weak solution of the closed loop equation, satisfying in addition
   \begin{equation}\label{5.6}\E \int_0^{\infty}   e^{- \lambda t}  |\underline{u}(X_t )|^2 dt < \infty.
   \end{equation}
\end{prop}
\begin{proof}
We start by constructing a canonical version of a cylindrical
Wiener process in $\Xi$. An explicit construction is needed to
clarify the application of an infinite-dimensional version of the
Girsanov theorem that we use below. We
 choose a larger Hilbert space $\Xi^{'}\supset \Xi$ in such a way that
 $\Xi$ is continuously and densely embedded in $\Xi^{'}$ with
 Hilbert-Schmidt inclusion operator $\mathcal{J}$. By $\Omega$ we denote the
 space $C([0,\infty[,\Xi^{'})$ of continuous functions $\omega: [0,\infty[ \rightarrow \Xi^{'}$
 endowed with the usual locally
 convex topology that makes $\Omega$ a Polish space, and by $\mathcal{B}$ its Borel
  $\sigma$-field. Since $\mathcal{J} \mathcal{J}^*$ has finite trace on $\Xi^{'}$, it is well known
that there exists a probability
 $\mathbb{P}$ on $\mathcal{B}$ such that the canonical processes
 $W^{'}_t(\omega):=\omega(t)$, $t \geq 0$, is a Wiener process with continuous paths in
 $\Xi^{'}$ satisfying $\mathbb{E}[\<W^{'}_t,\xi^{'}\>_{\Xi^{'}}\<W^{'}_s,\eta^{'}\>_{\Xi^{'}}]
 =\<\mathcal{J}\mathcal{J}^*\xi^{'},\eta^{'}\>_{\Xi^{'}}
 (t \wedge s) $
  for all $\xi^{'}, \eta^{'}\in \Xi^{'}$, $t,s \geq 0$.
This is called a $\mathcal{J} {\mathcal{J}}^*$-Wiener processes in
$\Xi^{'}$ in \cite{dz2}, to which we refer the reader for
preliminary material on Wiener processes on Hilbert spaces. Let us
denote by $\m G$ the $\mathbb{P}$-completion of $\mathcal{B}$
  and by $\m N$ the family of sets $A \in \m G$ with $\P(A)=0$. Let $\mathcal{B}_t=\sigma\{W^{'}_s: s\in [0,t]\}$ and $\mathcal{F}_t = \sigma(\mathcal{B}_t,\m N)$, $t\ge 0$,
  where as usual $\sigma(\cdot)$ denotes the $\sigma$-algebra in $\Omega$ generated by the indicated collection of sets or random variables. Thus $(\m F_t)_{t \geq 0}$ is the Brownian filtration of $W^{'}$.

  The $\Xi$-valued cylindrical Wiener process $\{W^{\xi}_t: t\geq 0, \xi\in \Xi\}$
  can now be defined as follows. For $\xi$ in the image of $\mathcal{J}^* \mathcal{J}$
  we take $\eta$ such that $\xi=\mathcal{J}^* \mathcal{J} \eta$ and define
  $W^{\xi}_s=\<W^{'}_s,\mathcal{J}\eta\>_{\Xi^{'}}$. Then we notice that
$\mathbb{E}|W^{\xi}_t|^2=t |\mathcal{J}\eta|_{\Xi^{'}}^2 = t
|\xi|_{\Xi}^2$, which shows that the mapping $\xi \rightarrow
W^{\xi}_s$, defined for $\xi \in \mathcal{J}^*
\mathcal{J}(\Xi)\subset \Xi$ with values in $L^2(\Omega,\m F,\P)$,
is an isometry for the norms of $\Xi$ and $L^2(\Omega,\m F,\P)$.
Consequently, noting that $\mathcal{J}^* \mathcal{J} (\Xi)$ is
dense in $\Xi$, it extends to an isometry $\xi \rightarrow
L^2(\omega,\m F,\P)$, still denoted $\xi \rightarrow W^{\xi}_s$.
An appropriate modification of $\{W^{\xi}_t: t\geq 0, \xi\in
\Xi\}$ gives the required cylindrical Wiener process. We note that
the Brownian filtration of $W$ coincides with $(\m F_t)_{t \geq
0}$.

Now let $X\in L^p_{ {\rm loc}}(\Omega,C(0,+\infty;H))$ be the mild
solution of
\begin{equation}\label{eqnocontr}
  \left\{
 \begin{array}{l}
 {\displaystyle dX_{\tau}=AX_{\tau}\;d\tau+b(X_{\tau})\;d\tau
 +\sigma \;dW_{\tau} } \\
    X_{0}=x
 \end{array}
 \right.
\end{equation}
If together with previous forward equation we also consider the
backward equation
\begin{equation}\label{back-system-contr}
Y_t -Y_T +\int_t^T Z_s dW_s = \int_t^T F(X_s,Y_s,Z_s)ds \qquad
0\leq t \leq T <\infty
\end{equation}
we know that there exists a unique solution $\{X_t^x,Y_t^x,Z_t^x,
t \geq 0\}$ forward-backward system
(\ref{eqnocontr})-(\ref{back-system-contr}) and by Proposition
\ref{main},
\begin{equation*}
    v(x) =Y_0^x.
    \label{relazionecontrollo1}
\end{equation*}
is the solution of the of the non linear stationary Kolmogorov
equation:
\begin{equation}
\call v(x) + F (x,v(x),\nabla v(x)\, \sigma)=0,\qquad x\in H.
\end{equation}
Moreover the following holds:
\begin{equation}\label{relazionecontrollo2}
Y_{\tau} (x)=v( X_\tau(x)), \qquad Z_\tau(x)=\nabla v(
X_\tau(x))\sigma
\end{equation}
We have \begin{equation}\label{5.7}\E \int_0^{\infty}e^{- (\lambda
+ \epsilon) t} |Z_t |^2 dt < \infty.\end{equation} and hence
                                  \begin{equation}\label{5.8}\E \int_0^T |Z_t |^2 dt < \infty.\end{equation}
By (\ref{r-gro}) we have \begin{equation}|r(X_t , \underline{u}(
X_t ))| \leq C(1 + |\underline{u}( X_t )|),\end{equation} and by
(\ref{5.3}),
\begin{equation}\label{5.9}|\underline{u}( X_t )| = |\gamma(X_t ,\nabla v(
X_t(x))\sigma)| \leq C(1  + |\nabla v( X_t(x))\sigma|) = C(1  +
|Z_t |).
\end{equation} Let us define $\forall
T>0$
 \begin{equation}\label{defdirho-bis}
M_T=\exp\left( \int_0^T\langle
 r(X_s,\underline{u}(X_s),dW_s\rangle_{\Xi} -\frac{1}{2}\int_0^T
|r(X_s,\underline{u}(X_s) |^2_{\Xi}\; ds \right).
\end{equation}
Now, arguing exactly as in the proof of Proposition 5.2 in
\cite{FHT}, we can prove that $\E M_T = 1$, and $M$ is a
$\P$-martingale. Hence there exists a probability
$\widehat{\mathbb{P}}_T$ on $\mathcal{F}_T$ admitting $M_T$ as a
density with respect to $\mathbb{P}$, and by the Girsanov Theorem
we can conclude that $\{\widehat{W}_t, t \in [0,T] \}$ is a Wiener
process with respect to $\P$ and $(\m F_t )$.
Since $\Xi^{'}$ is a Polish space and $\widehat{\mathbb{P}}_{T+h}$
coincide with $\widehat{\mathbb{P}}_{T}$ on $\mathcal{B}_T$, $T,h
\geq 0$, by known results (see \cite{RY}, Chapter VIII, \S 1,
Proposition (1.13)) there exists a probability
$\widehat{\mathbb{P}}$ on $\mathcal{B}$ such that the restriction
on $\mathcal{B}_T$ of $\widehat{\mathbb{P}}_T$ and that of
$\widehat{\mathbb{P}}$ coincide, $T\geq 0$. Let $\widehat{\m G}$
be the $\widehat{\mathbb{P}}$-completion of $\mathcal{B}$ and
$\widehat{\mathcal{F}}_T$ be the $\widehat{\mathbb{P}}$-completion
of $\mathcal{B}_T$. Moreover, since  for all $T>0$,
$\{\widehat{W}_{t}: t\in [0,T]\}$ is a $\Xi$-valued cylindrical
Wiener process under $\widehat{\mathbb{P}}_T$ and the restriction
of $\widehat{\mathbb{P}}_T$ and of $\widehat{\mathbb{P}}$ coincide
on $\mathcal{B}_T$ modifying $\{\widehat{W}_{t} : t\geq 0\}$ in a
suitable way on a $\widehat{\mathbb{P}}$-null probability set we
can conclude that $(\Omega,\widehat{{\m
G}},\{\widehat{\mathcal{F}}_{t}, t\geq 0\},\widehat{{\P}},
\{\widehat{{W}}_{t}, t\geq 0\},\gamma(X,
  \nabla v(X)\sigma(X)) )$ is an admissible control system. The above construction immediately ensures that,
   if we choose such an admissible control system, then (\ref{closedloop}) is satisfied.
  Indeed if we rewrite (\ref{eqnocontr}) in terms of $\{\widehat{W}_{t} : t\geq 0\}$ we get
 $$
  \left\{
 \begin{array}{l}
 {\displaystyle dX_{\tau}=AX_{\tau}\;d\tau+b(X_{\tau})\;d\tau
 +\sigma \;[r(X_{\tau},\underline{u}(X_\tau) )d\tau+d\widehat{W}_{\tau}] } \\
    X_{0}=x.
 \end{array}
 \right.$$

It remains to prove (\ref{5.6}). We define stopping times
 $$\sigma_n = \inf\bigg\{t\geq 0 :      \int_0^t  e^{- \lambda t} |Z_s |^2 ds \geq n \bigg\},$$
with the convention that $\sigma_n = \infty$ if the indicated set
is empty. By (\ref{5.7}) for $\P$-a.s. $\omega \in \Omega$ there
exists an integer $N(\omega)$ depending on $\omega$ such that
$\sigma_n (\omega) = \infty$ for $n \geq N (\omega)$. Applying the
Ito formula to $e^{- \lambda t}Y_t$, with respect to W , we obtain
$$ \begin{array}{l}\dis
e^{- \lambda \sigma_n} Y_{\sigma_n} = Y_0  -\int^{\sigma_n}_0 e^{-
\lambda s }Z_s \;d{W}_s+
 \\\dis\qquad
 + \int^{\sigma_n}_0 e^{- \lambda s }\left[-F(X_s,Y_s,
  Z_s) -  \lambda  Y_s (x)ds+
Z_sr(X_s,\underline{u}( X_s )) \right] \;ds.
 \end{array}
$$
from which we deduce that
 $$  \begin{array}{l}\dis
\E e^{- \lambda \sigma_n} Y_{\sigma_n}  +\E\int^{\sigma_n}_0  e^{-
\lambda s }g(X_s , \underline{u}( X_s )) ds  = Y_0  +
 \\\dis\qquad
 + \E\int^{\sigma_n}_0 e^{- \lambda s }\left[-F(X_s,Y_s,
  Z_s) - \lambda  Y_s ds+
Z_sr(X_s,\underline{u}( X_s )) + g(X_s , \underline{u}( X_s
))\right] ds= Y_0 .
 \end{array}
$$
with the last equality coming from the definition of
$\underline{u}$. Recalling that $Y$ is bounded, it follows that
$$       E\int^{\sigma_n}_0 e^{-\lambda s}g(X_s , \underline{u}( X_s )) ds \leq C                      $$
for some constant $C$ independent of $n$. By (\ref{2.10}) and by
sending $n$ to infinity, we finally prove (\ref{5.6}).

\end{proof}


\begin{thebibliography}{10}

\bibitem{AP} A. Ambrosetti, G. Prodi.
\newblock{\em A primer of nonlinear analysis},
\newblock Cambridge Studies in Advanced Mathematics, 34,
Cambridge University Press, 1995.

\bibitem{BC} Ph. Briand, F. Confortola.
\newblock BSDEs with stochastic Lipschitz condition and quadratic PDEs in Hilbert spaces.
\newblock{\em Stochastic Processes and their Applications.} To appear. 

\bibitem{BH98}  P. Briand, Y. Hu.
\newblock Stability of BSDEs with random terminal time and
homogenization of semilinear elliptic PDEs.
\newblock{\em J. Funct. Anal.} {\bf 155} (1998), 455-494.


\bibitem{BuPe}  R. Buckdahn, S. Peng.
\newblock Stationary backward stochastic differential equations and
associated partial differential equations.
\newblock {\em Probab. Theory Related
Fields} {\bf 115} (1999), 383-399.

\bibitem{Ce}  S. Cerrai,
\newblock{\em Second order PDE's in finite and infinite dimensions. A
probabilistic approach}.
\newblock Lecture Notes in Mathematics
{\bf 1762}, Springer, Berlin, 2001.

\bibitem{DaPa97}  R. W. R. Darling, E. Pardoux.
 \newblock Backwards SDE with random terminal time and applications to
semilinear elliptic PDE,
\newblock {\it Ann. Probab.} {\bf 25} (1997),
1135-1159.
\bibitem{ElK97}
N.~El~Karoui.
\newblock Backward stochastic differential equations: a general introduction.
\newblock In {\em Backward stochastic differential equations (Paris,
  1995--1996)}, volume 364 of {\em Pitman Res. Notes Math. Ser.}, pages 7--26.
  Longman, Harlow, 1997.
\bibitem{dz2} G. Da Prato and J. Zabczyk.
\newblock {\em Stochastic Equations in Infinite
Dimensions.} {Cambridge University Press, Cambridge, 1992}

\bibitem{dz3} G. Da Prato and J. Zabczyk.
\newblock {\em Second Order Partial Differential Equations
in Hilbert Spaces}.
\newblock Cambridge University Press, Cambridge, 2002.

\bibitem{fede} F. Masiero.
\newblock Infinite horizon stochastic optimal control problems with degenerate noise and
elliptic equations in Hilbert spaces. Preprint, Politecnico di
Milano, 2004 (submitted).

\bibitem{FHT} M. Fuhrman,Y.~Hu and G. Tessitore.
\newblock On a class of stochastic optimal control problems related to BSDEs with quadratic growth.
\newblock{\em SIAM J. Control Optim.} {\bf 45}  (2006),  no. 4, 1279--1296.

\bibitem{FT1}  M. Fuhrman and G. Tessitore.
\newblock Nonlinear Kolmororov equations
in infinite dimensional spaces: the backward stochastic
differential equations approach and applications to optimal
control.
\newblock{\em Ann. Probab.} {\bf 30} (2002), 1397-1465.


\bibitem{FT3}  M. Fuhrman and G. Tessitore,
\newblock Infinite horizon backward stochastic differential equations
 and elliptic equations in Hilbert spaces.
\newblock{\em Ann. Probab.} {\bf 30} (2004), 607-660.

\bibitem{GoRo} F. Gozzi, E. Rouy.
\newblock Regular solutions of second-order stationary
              {H}amilton-{J}acobi equations.
\newblock{\em J. Differential Equations} {\bf 130} (1996), 201-234.

\bibitem{Kaz94} N.~Kazamaki.
\newblock {\em Continuous exponential martingales and {BMO}}, volume 1579 of
  {\em Lecture Notes in Mathematics}.
\newblock Springer-Verlag, Berlin, 1994.

\bibitem{HIM05}
Y.~Hu, P.~Imkeller, and M.~M{\"u}ller.
\newblock Utility maximization in incomplete markets.
\newblock {\em Ann. Appl. Probab.}, 15(3):1691--1712, 2005.

\bibitem{HT}
Y. Hu and G. Tessitore.
\newblock BSDE on an infinite horizon and elliptic PDEs in infinite dimension.
\newblock {\em NoDEA Nonlinear Differential Equations Appl.} To appear.

\bibitem{Kob00}
M.~Kobylanski.
\newblock Backward stochastic differential equations and partial differential
  equations with quadratic growth.
\newblock {\em Ann. Probab.}, 28(2):558--602, 2000.

\bibitem{LsM98}
J.-P. Lepeltier and J.~San~Martin.
\newblock Existence for {B}{S}{D}{E} with superlinear-quadratic coefficient.
\newblock {\em Stochastic Stochastics Rep.}, 63(3-4):227--240, 1998.




\bibitem{Pa98} E. Pardoux.
\newblock Backward stochastic differential equations and viscosity
solutions of systems of semilinear parabolic and elliptic PDEs of
second order, in: {\em Stochastic Analysis and related topics},
the Geilo workshop 1996, eds. L. Decreusefond, J. Gjerde, B
{\O}ksendal, A.S. \"Ust\"unel, 79-127, Progress in Probability
{\bf 42}, Birkh\"auser, Boston, 1998.

\bibitem{Pa99} {\'E}.~Pardoux.
\newblock BSDEs, weak convergence and homogenization
of semilinear PDEs.
\newblock {\em Nonlinear analysis, differential equations and control (Montreal,
QC, 1998), 503--549, NATO Sci. Ser. C Math. Phys. Sci., 528, }
\newblock Kluwer Acad. Publ., Dordrecht, 1999.

\bibitem{RY} D. Revuz and M. Yor.
\newblock{\em Continuous martingales and Brownian motion.}
\newblock Grundlehren der Mathematischen Wissenschaften [Fundamental Principles of Mathematical Sciences], 293. Springer-Verlag, Berlin, (1999).

\bibitem{Ro}  M. Royer.
\newblock BSDEs with a random terminal time driven by
 a monotone generator and their links with PDEs.
\newblock {\em Stochastics Stochastics Rep.}
{\bf 76} (2004) 281-307.
\end{thebibliography}
\end{document}